\documentclass{amsart}
\usepackage[a4paper,twoside,inner=2cm,outer=2cm]{geometry}
\usepackage{amssymb,amsmath,amscd}
\usepackage{tikz}
\usetikzlibrary{matrix,arrows}
\usepackage{graphicx}  
\usepackage{float}  
\usepackage{multirow}  
\usepackage{framed}
\usepackage{booktabs}  
\usepackage{subcaption}  
\usepackage[hidelinks,colorlinks]{hyperref}
\usepackage{lineno}

\usepackage[shortlabels]{enumitem}

\theoremstyle{definition}
\newtheorem{dfn}{Definition}[section]

\theoremstyle{remark}
\newtheorem{rmk}[dfn]{Remark}
\newtheorem{exm}[dfn]{Example}
\theoremstyle{plain}
\newtheorem{thm}[dfn]{Theorem}
\newtheorem*{thm*}{Theorem}
\newtheorem{prop}[dfn]{Proposition}
\newtheorem{lem}[dfn]{Lemma}
\newtheorem{cor}[dfn]{Corollary}

\newcommand*{\Cdot}{\raisebox{-1ex}{\scalebox{3}{$\cdot$}}}
\renewcommand{\leq}{\leqslant}
\renewcommand{\geq}{\geqslant}
\renewcommand{\setminus}{\smallsetminus}

\newcommand{\Z}{\mathbb{Z}}
\newcommand{\Q}{\mathbb{Q}}
\newcommand{\R}{\mathbb{R}}
\newcommand{\C}{\mathbb{C}}
\newcommand{\F}{\mathbb{F}}
\newcommand{\K}{\mathbb{K}}

\newcommand{\D}{\mathbb{D}}

\newcommand{\longhookrightarrow}{\lhook\joinrel\longrightarrow}

\title{Equivariant cohomology rings of the real flag manifolds}



\author[He]{\bfseries Chen He}

\address{
	Yau Mathematical Sciences Center \\ 
	Tsinghua University   \\ 
	Beijing\\
	P.R.\,China}

\email{che@math.tsinghua.edu.cn}

\begin{document}

\vspace{18mm} \setcounter{page}{1} \thispagestyle{empty}

\begin{abstract}
	We give Leray-Borel-type descriptions for the mod-$2$ and the rational equivariant cohomology rings of the real and the oriented flag manifolds under the canonical torus or 2-torus actions.
\end{abstract}

\maketitle

\tableofcontents

\section{Introduction}
\vskip 15pt

Let $G$ be a compact connected Lie group and $H$ be a closed connected subgroup. A general method of computing the cohomology of the homogeneous space $G/H$ in rational coefficients was given by H.\,Cartan\,\cite{Car50} as the cohomology of a cochain complex $\big(H^*(G,\Q)\otimes H^*(BH,\Q),\,\delta\big)$, where $BH$ is the classifying space of $H$ and $\delta$ is a differential. Though theoretically feasible, computing the cohomology of homogeneous spaces is not always straightforward, especially when accounting for the ring structures and in coefficient rings of $\Z$ or $\F_p$ the finite field of a prime $p$.

Nevertheless, several major types of homogeneous spaces have very explicit descriptions of their cohomology. When $H$ is the identity group, the homogeneous space $G/H$ is the Lie group $G$, whose rational homology and cohomology have been shown to be exterior algebras by Hopf\,\cite{Ho41}. When $G$ and its subgroup $H$ are of the same rank (we will call this the \textit{equal-rank condition}), the cohomology ring of $G/H$ can be described as the quotient of Weyl group invariant polynomial rings, proved by Leray\,\cite{Le50} for rational coefficients, and by Borel\,\cite{Bo53a} for other appropriate coefficient rings. Moreover, if $G$ and $H$ contain the same maximal torus $T$, then the canonical $T$-action on $G/H$ from the left by multiplications are equivariantly formal. Descriptions of various $T$-equivariant cohomology theories for such $G/H$ can be found for example in  \cite{KV93,Br97,GHZ06,KiKr13,CZZ15}, and those when $G$ is a Kac-Moody group can be found for example in \cite{Ar89,KK90,HHH05}. 

An important family of homogeneous spaces are the flag manifolds in fields $\D=\R,\C,\mathbb{H}$ for a sequence of positive integers $n_1,\ldots,n_k$ and $n=n_1+\cdots+n_k$: the space $Fl(n_1,\ldots,n_k,\D^n)={G(n)}/\big(G(n_1)\times\cdots\times G(n_k)\big)$, where $G=O,\,U,\,Sp$, consists of the $\D$-orthogonal decompositions $W_1\oplus\cdots\oplus W_k=\D^n$ with $\dim_\D W_i = n_i$, or equivalently consists of the flags $V_1 \subset V_2 \cdots \subset V_{k-1}  \subset V_k =\D^n$ of dimensions $n_1<n_1+n_2<\cdots<n_1+\cdots+n_{k-1}<n$. The oriented flag manifold will be denoted by $Fl^o(n_1,\ldots,n_k,\R^n)={SO(n)}/\big(SO(n_1)\times\cdots\times SO(n_k)\big)$. Among these flag manifolds, the complex or quaternionic ones satisfy the equal-rank condition. Their integral cohomology rings in terms of characteristic classes were given by Borel\,\cite{Bo53a}. However, many real or oriented flag manifolds, for example, the odd-dimensional real or oriented Grassmannians and the complete real or oriented flag manifolds do not satisfy the equal-rank condition. The rational cohomology rings of some cases of the oriented flag manifolds have been given in the book of Greub, Halperin and Vanstone\,\cite{GHV76}. 

Another important family of homogeneous spaces are the Stiefel manifolds of $\D$-orthonormal $k$-frames in $\D^n$: $V_k(\D^n)=G(n)/G(n-k)$ whose cohomology rings were given by Borel\,\cite{Bo53a} and Miller\,\cite{Mi53}. But most of the Stiefel manifolds do not satisfy the equal-rank condition. 

In order to give a unified treatment of the above concrete examples regardless of the equal-rank condition, we will consider the homogeneous space ${G(n)}/{\big(G(n_1)\times\cdots\times G(n_k)\big)}$ where the symbol $G=O,\,U,\,Sp,\,SO$ and $n_i\geq 1$ but only requiring $n\geq n_1+\cdots+n_k$. We will still call them \textit{flag manifolds} and denote them by $Fl(n_1,n_2,\ldots,n_k,\D^n)$ for $\D=\R,\C,\mathbb{H}$ or $Fl^o(n_1,n_2,\ldots,n_k,\R^n)$. To distinguish, the case of $n= n_1+\cdots+n_k$ will be called the \textit{usual flag manifolds}. 

In this paper, we will consider the canonical actions of tori or $2$-tori on the real and the oriented flag manifolds and give Leray-Borel-type descriptions of their equivariant cohomology rings in terms of equivariant characteristic classes together with some additional generators.

\textit{Acknowledgement}. The author is very grateful to Jeff Carlson, Oliver Goertsches, Liviu Mare and Xu'an Zhao for many valuable discussions. The author thanks the Special Support of the China Postdoctoral Science Foundation.

\vskip 20pt
\section{Preliminaries}
\vskip 15pt

In this section, we consider some simple constructions on the real or oriented flag manifolds and review some basic notions of the equivariant cohomology. 

\subsection{Real flag manifolds}
Given a sequence of positive integers $n_1,\ldots,n_{k}$ and an integer $n\geq n_1+\cdots+n_{k}$, we consider the real flag manifold  $Fl(n_1,\ldots,n_k,\R^n)\triangleq O(n)/\big(O(n_1)\times \cdots \times O(n_k)\big)$ and the oriented flag manifold $Fl^o(n_1,\ldots,n_k,\R^n)\triangleq SO(n)/\big(SO(n_1)\times \cdots \times SO(n_k)\big)$. 
\subsubsection{A bundle structure}\label{subsubsec:generalFlag}
Set $N=n_1+\cdots+n_{k}$. The real flag manifold has a natural bundle structure with the fibre a usual real flag manifold and the base a real Stiefel manifold:
\begin{align*}
\frac{O(N)}{O(n_1)\times \cdots \times O(n_k)}
\overset{\iota}{\longhookrightarrow}
\frac{O(n)}{O(n_1)\times \cdots \times O(n_k)}
\overset{\pi}{\longrightarrow}
\frac{O(n)}{O(N)}.
\end{align*}
It enables us to write an element of $Fl(n_1,\ldots,n_k,\R^n)$ as a combination of an $(n-N)$-frame and an orthonormal decomposition:
\begin{gather*}
(v_1,\ldots,v_{n-N};W_1,\ldots,W_{k}) \\
 \mbox{s.t. } v_i\in \R^n, \langle v_i,v_j\rangle =\delta_{i,j}\, \forall 1\leq i, j \leq n-N;\, W_i\subset \R^N,\langle W_i,W_j\rangle=0,\, \forall 1\leq i\neq j \leq k
\end{gather*}
where $W_i$ is of dimension $n_i$ and $\oplus_{i=1}^k W_i=(\mathrm{Span}_\R\{v_1,\ldots,v_{n-N}\})^{\perp}$. Such an element will be denoted by $(v_\bullet;W_\bullet)$. The projection $\pi:Fl(n_1,\ldots,n_k,\R^n)\rightarrow V_{n-N}(\R^n)$ sends $(v_\bullet;W_\bullet)$ to $v_\bullet$. Write $\R^n=\R^{N}\oplus \R^{n-N}$, and let $v^0_\bullet$ be the standard orthonormal basis of $\R^{n-N}$, an inclusion $\iota:Fl(n_1,\ldots,n_k,\R^N)\rightarrow Fl(n_1,\ldots,n_k,\R^n)$ can be chosen to send $W_\bullet$ to $(v^0_\bullet;W_\bullet)$.

The oriented flag manifold $Fl^o(n_1,\ldots,n_k,\R^n)= SO(n)/\big(SO(n_1)\times \cdots \times SO(n_k)\big)$ has a similar bundle structure. Notice that $Fl^o(n_1,\ldots,n_k,\R^n)=SO(n)/\big(SO(n_1)\times \cdots \times SO(n_k)\times SO(1)^{n-N}\big)$ is a usual oriented flag manifold in a degenerate sense.

\subsubsection{Universal vector bundles}\label{subsubsec:univBundle}
There are natural vector bundles on a real flag manifold $Fl(n_1,\ldots,n_k,\R^n)$:
\begin{align*}
\mathcal{L}_i &\triangleq \big\{ (v_\bullet;W_\bullet;v) \in Fl(n_1,\ldots,n_k,\R^n)\times \R^n \mid v=zv_i \mbox{ for some } z\in \R \big\} & & 1\leq i\leq n-N\\
\mathcal{W}_j &\triangleq \big\{(v_\bullet;W_\bullet;v)\in Fl(n_1,\ldots,n_k,\R^n)\times \R^n \mid v\in W_j\big\}& & 1\leq j\leq k.
\end{align*} 

Each line bundle $\mathcal{L}_i$ is a trivial bundle. Each vector bundle $\mathcal{W}_j$ is of rank $n_j$, and will be called the \textit{$j$-th universal vector bundle} over $Fl(n_1,\ldots,n_k,\R^n)$.

Consider the \textit{full universal vector bundle} 
\begin{align}\label{eq:fullDecom}
\mathcal{V}\triangleq\oplus_{j=1}^k \mathcal{W}_j.
\end{align}
By definition, the right hand side satisfies $\oplus_{j=1}^k \mathcal{W}_j = \big(\oplus_{i=1}^{n-N} \mathcal{L}_i\big)^\perp$, hence we have the \textit{complete decomposition}:
\begin{align}\label{eq:completeDecom}
Fl(n_1,\ldots,n_k,\R^n)\times \R^n = \big(\oplus_{i=1}^{n-N} \mathcal{L}_i\big) \oplus \mathcal{V}= \big(\oplus_{i=1}^{n-N} \mathcal{L}_i\big) \oplus \big(\oplus_{j=1}^k \mathcal{W}_j \big).
\end{align}

Let ${w}^{(j)}=1+{w}_1^{(j)}+\cdots+{w}_{n_j}^{(j)}$ and ${p}^{(j)}=1+{p}_1^{(j)}+\cdots+{p}_{[n_j/2]}^{(j)}$ be the total Stiefel-Whitney and Pontryagin classes of the $j$-th universal vector bundle $\mathcal{W}_j$ over the real flag bundle. The Whitney product formula of total characteristic classes, applied to the complete decomposition\,\eqref{eq:completeDecom}, gives the relations
\begin{gather}\label{eq:Whitney}
\prod_{j=1}^{k}{w}^{(j)}=1, \qquad \qquad \qquad  \prod_{j=1}^{k}{p}^{(j)}=1.
\end{gather}

Over the oriented flag manifold $Fl^o(n_1,\ldots,n_k,\R^n)$, the universal vector bundles are canonically oriented. Let ${e}^{(j)}$ be the Euler class of $\mathcal{W}_j$. The square formula of Euler class gives
\begin{align}\label{eq:sqEuler}
(e^{(j)})^2=p^{(j)}_{[n_j/2]}
\end{align}
where $p^{(j)}_{[n_j/2]}$ is the top Pontryagin class of $\mathcal{W}_j$.
The Whitney product formula of Euler classes, applied to decomposition~\eqref{eq:fullDecom}, gives the relation
\begin{align}\label{eq:WhitneyEuler}
\prod_{j=1}^{k}e^{(j)}=e(\mathcal{V})
\end{align}
where $e(\mathcal{V})$, the Euler class of the full universal bundle $\mathcal{V}$, will be specified shortly in later sections.

From the inclusion of groups $O(k-1)\subset O(k)\subset O(n)$, we have the bundle of Stiefel manifolds
\[
S^{k-1}=O(k)/O(k-1)\longrightarrow O(n)/O(k-1)\overset{\pi}{\longrightarrow} O(n)/O(k).
\]

\begin{lem}\label{lem:BundleTower}
	The map $\pi: O(n)/O(k-1)\rightarrow O(n)/O(k)$ is the projection $V_{n-k+1}(\R^{n})\rightarrow V_{n-k}(\R^n)$ sending an $(n-k+1)$-frame $(v_1,\ldots,v_{n-k},v_{n-k+1})$ to the $(n-k)$-frame $(v_1,\ldots,v_{n-k})$. Moreover,
	\begin{enumerate}[label=(\alph*)]
		\item the bundle $S^{k-1}\hookrightarrow V_{n-k+1}(\R^{n})\rightarrow V_{n-k}(\R^n)$ is the unit sphere bundle of the rank-$k$ full universal vector bundle $\mathcal{V}\big(V_{n-k}(\R^{n})\big)\rightarrow V_{n-k}(\R^{n})$,
		\item the kernel of the tangent map $T\,\pi:T\,V_{n-k+1}(\R^{n})\rightarrow T\,V_{n-k}(\R^n)$ is the vector bundle over $V_{n-k+1}(\R^{n})$ with fibres being the tangent spaces along the sphere fibres of $\pi$ and we have $\ker\,(T\,\pi)\cong \mathcal{V}\big(V_{n-k+1}(\R^{n})\big)$.
	\end{enumerate}
\end{lem}
\begin{proof}
	The first statement is straightforward.
	\begin{enumerate}[(a)]
		\item Take any orthonormal $(n-k)$-frame $(v_1,\ldots,v_{n-k})\in V_{n-k}(\R^n)$, it extends to an orthonormal $(n-k+1)$-frame $(v_1,\ldots,v_{n-k},v_{n-k+1})\in V_{n-k+1}(\R^n)$ iff the unit vector $v_{n-k+1}$ is in $\big(\mathrm{span}_\R\,\{v_1,\ldots,v_{n-k}\}\big)^\perp$ which is a fibre of the full universal vector bundle $\mathcal{V}\big(V_{k}(\R^{n})\big)$. 
		\item The identification of $\ker\,(T\,\pi)$ as the fiber-wise tangent bundle is straightforward. For the second identification, take any $(v_1,\ldots,v_{n-k})\in V_{n-k}(\R^{n})$. By (a), the sphere fibre of $\pi$ at $(v_1,\ldots,v_{n-k})$ is  $$\pi^{-1}(v_1,\ldots,v_{n-k})=\big\{(v_1,\ldots,v_{n-k},v_{n-k+1}) \mid v_{n-k+1} \in (\mathrm{span}_\R\,\{v_1,\ldots,v_{n-k}\})^\perp,  \|v_{n-k+1}\|=1 \big\}. $$ 
		The fiber of $\ker\,(T\,\pi)$ at $(v_1,\ldots,v_{n-k},v_{n-k+1})$, as the tangent space of $\pi^{-1}(v_1,\ldots,v_{n-k})$ at $v_{n-k+1}$, is 
		$$\big(\mathrm{span}_\R\,\{v_{n-k+1}\}\big)^\perp\cap \big(\mathrm{span}_\R\,\{v_1,\ldots,v_{n-k}\}\big)^\perp=\big(\mathrm{span}_\R\,\{v_1,\ldots,v_{n-k},v_{n-k+1}\}\big)^\perp$$ 
		which is exactly the fibre of $\mathcal{V}\big(V_{n-k+1}(\R^{n})\big)$ at $(v_1,\ldots,v_{n-k},v_{n-k+1})$.
	\end{enumerate}
\end{proof}

\subsubsection{Canonical torus actions} \label{subsubsec:canTorusAction}
Suppose we have a sequence $(n_1,\ldots,n_k)=(2m_1+1,\ldots,2m_{l}+1,2m_{l+1},\ldots,2m_k)$ with $l$ odd integers and $k-l$ even integers. The product $O(2m_1+1)\times \cdots \times O(2m_{l}+1)\times O(2m_{l+1})\times \cdots \times O(2m_k)$ has a maximal torus $T^M$ of rank $M=m_1+\cdots +m_k$. The torus acts by left multiplications on the real and oriented flag manifolds. We use the notation of multinomial coefficient $\binom{m}{m_1,\ldots,m_k}\triangleq m!/(m_1!\cdots m_k!)$.

\begin{lem}\label{lem:fixedPointOfRealGeneralFlag}
	The fixed-point sets of the $T^M$-actions on the real or oriented flag manifolds $O(n)/\big(O(n_1)\times \cdots \times O(n_k)\big),\,SO(n)/\big(SO(n_1)\times \cdots \times SO(n_k)\big)$ are respectively a disjoint union of $\binom{M}{m_1,\ldots,m_k}$ copies of
	\[
	Fl(\underbrace{1,\ldots,1}_{l \text{ items}},\R^{n-2M})=\frac{O(n-2M)}{O(1)^l}  \qquad \mbox{or} \qquad 
	Fl^o(\underbrace{1,\ldots,1}_{l \text{ items}},\R^{n-2M})=\frac{SO(n-2M)}{SO(1)^l}=SO(n-2M).
	\]
\end{lem} 
\begin{proof}
	We work on the real case. The $T^M$-action on $Fl(2m_1+1,\ldots,2m_{l}+1,2m_{l+1},\ldots,2m_k;\R^n)$ is induced from the real representation of $T^M$ on $\R^n=\big\{(\R^2)^{m_1}\oplus \R\big\}\oplus\cdots\oplus\big\{(\R^2)^{m_l}\oplus \R\big\}\oplus (\R^2)^{m_{l+1}}\oplus \cdots \oplus(\R^2)^{m_{k}} \oplus \R^{n-2M-l}$
	decomposed into irreducible parts such that each $(\R^2)^{m_i}$ is acted by a $T^{m_i}$-factor of $T^M$, but the individual $\R$-summands and the $\R^{n-2M-l}$-summand are trivial representations. Rearranging the summands, we can write $\R^n= (\R^2)^M \oplus \R^{n-2M}$ as a direct sum of the standard real representation of $T^M$ and a trivial representation.
	
	Let $(v_1,\ldots,v_{n-N};W_1,\ldots,W_{k})\in Fl(n_1,\ldots,n_k;\R^n)$ be fixed by the $T^M$-action, where $N=n_1+\cdots n_k = 2M+l$. Firstly, we have every $v_i\in (\R^n)^{T^M}=\R^{n-2M}$, hence $(v_1,\ldots,v_{n-2M-l})\in V_{n-2M-l}(\R^{n-2M})$. Secondly, let $U=\big(\mathrm{Span}_\R\{v_1,\ldots,v_{n-2M-l}\}\big)^{\perp}\cap \R^{n-2M}$, then
	$\big(\mathrm{Span}_\R\{v_1,\ldots,v_{n-2M-l}\}\big)^{\perp}= (\R^2)^{M}\oplus U$ as a direct sum of a standard $T^M$-representation and a trivial representation.
	Since $\oplus_{i=1}^k W_i=\big(\mathrm{Span}_\R\{v_1,\ldots,v_{n-2M-l}\}\big)^{\perp}$, there is a partition $\{1,2,\ldots,{M}\}=S_1\cup S_2\cup\cdots\cup S_k$ into $k$ subsets of cardinalities $m_1,\ldots,m_k$, such that for $1\leq i\leq l$ we have $W_i=(\oplus_{S_i}\R^2) \oplus L_i$ with a line $L_i \subset U$ and $\oplus_{i=1}^{l}L_i = U$, and for $l+1\leq i\leq k$ we have $W_i=\oplus_{S_i}\R^2$. The number of such partitions is exactly $\binom{M}{m_1,\ldots,m_k}$. For every such partition, we have
	\[
	(v_1,\ldots,v_{n-2M-l};L_{1},\ldots,L_l)\in Fl(\underbrace{1,\ldots,1}_{l \text{ items}},\R^{n-2M})=\frac{O(n-2M)}{O(1)^l}.
	\]
    Conversely, every such $(v_\bullet;L_\bullet)$ together with an  $(n_1,\ldots,n_k)$-partition of $\{1,\ldots,M\}$ forms a $T^M$-fixed point.
\end{proof}

Let $N=n_1+\cdots+n_{k}$, the product $O(n_1)\times \cdots \times O(n_k)$ has the maximal $2$-torus $O(1)^N=\Z_2^N$. Similarly, we consider the action of $\Z_2^N$ by left multiplications on the real flag manifold.

\begin{lem}\label{lem:fixedPointOfRealGeneralFlag2}
	The fixed-point set of the $\Z_2^N$-action on the real flag manifold $O(n)/\big(O(n_1)\times \cdots \times O(n_k)\big)$ is a disjoint union of $\binom{n}{n_1,\ldots,n_k}$ copies of $O(n-N)$, i.e. $2\binom{n}{n_1,\ldots,n_k}$ copies of $SO(n-N)$.
\end{lem} 

\subsection{Equivariant cohomology}\label{subsec:equivCohom}
Let $G\curvearrowright X$ be a compact Lie group action on a compact manifold. 
\begin{dfn}
	The \textit{$G$-equivariant cohomology} of $X$ in coefficient ring $\K$ is $H^*_G(X,\K)\triangleq H^*(EG\times_G X,\K)$ where $EG$ is the universal $G$-bundle and the \textit{Borel construction} $EG\times_G X$ is the product $EG\times X$ taken quotient by the diagonal $G$-action.
\end{dfn}

By definition, for the trivial $G$-action on a single point $pt$, we have $H^*_G(pt,\K)=H^*(BG,\K)$ where $BG=EG/G$ is the classifying space of $G$. If $G\curvearrowright X$ is trivial, then $H^*_G(X,\K)\cong H^*_G(pt,\K)\otimes_\K H^*(X,\K)$. If $G\curvearrowright X$ is free, then $H^*_G(X,\K)\cong H^*(X/G,\K)$.

The fibration $X\overset{\iota}{\rightarrow} EG\times_G X \overset{\pi}{\rightarrow} EG/G=BG$ induces natural pullback homomorphisms of cohomology rings. The homomorphism $\pi^*:H^*(BG,\K){\rightarrow} H_G^*(X,\K)$ makes $H_G^*(X,\K)$ an $H^*(BG,\K)$-algebra. The homomorphism $\iota^*:H_G^*(X,\K){\rightarrow} H^*(X,\K)$ relates equivariant cohomology with ordinary cohomology. Given $\omega\in H^*(X,\K)$, if there exists $\tilde{\omega}\in H_G^*(X,\K)$ with $\omega =\iota^*(\tilde{\omega})$, then we say $\tilde{\omega}$ is an \textit{equivariant lift} of $\omega$.

\subsubsection{Equivariant formality}
The fibration $X\rightarrow EG\times_G X \rightarrow EG/G=BG$ produces a Leray-Serre spectral sequence with $E_2 =  H^*(BG,\K)\otimes_\K H^*(X,\K)$ that converges to $H^*(EG\times_G X,\K)=H^*_G(X,\K)$. 

\begin{dfn}\label{dfn:formality}
	An action $G\curvearrowright X$ is \textit{equivariantly formal} in coefficient $\K$ if the Leray-Serre spectral sequence of $X\rightarrow EG\times_G X \rightarrow BG$ collapses at $E_2 =  H^*(BG,\K)\otimes_\K H^*(X,\K)$ with an $H^*(BG,\K)$-module isomorphism
	\[
	H_G^*(X,\K)\cong H^*(BG,\K)\otimes_\K H^*(X,\K).
	\]
\end{dfn}

\begin{rmk}\label{rmk:formality}
	If both $H^*(X,\K),\,H^*(BG,\K)$ only have even-degree elements, then any action $G\curvearrowright X$ is equivariantly formal in coefficient $\K$. 
\end{rmk}

In the equivariantly formal situation, the ordinary cohomology can be recovered from the equivariant cohomology by the formula $H^*(X,\K)\cong H_G^*(X,\K)\otimes_{H^*(BG,\K)} \K$.

When the acting group is a torus or a $p$-torus, there is a numeric criterion for the equivariant formality.

\begin{thm}[Cohomology inequality and equivariant formality, see \cite{AP93} p.\,210, Thm\,3.10.4]\label{thm:CohomIneq}
	Let $G$ be a torus $(S^1)^n$ and set $\K=\Q$ or $\F_p$ of any prime $p$, or let $G$ be a discrete $p$-torus $\Z_p^n$ of a prime $p$ and set $\K=\F_p$. Suppose $G$ acts on $X$ and denote the fixed-point set by $X^G$, then the total Betti numbers in $\K$ coefficients satisfy the inequality
	$\sum_i \dim_\K\,H^i(X^G,\K) \leq \sum_i \dim_\K\,H^i(X,\K)$
	where the equality holds if and only if the action $G\curvearrowright X$ is equivariantly formal in coefficient $\K$.
\end{thm}

For the action of a torus $T$, we can use the $\Z$-ranks to replace the $\Q$-dimensions and get the inequality regarding the total Betti numbers of the free parts 
$
\sum_i \mathrm{rank}_\Z\,H^i(X^T,\Z) \leq \sum_i \mathrm{rank}_\Z\,H^i(X,\Z).
$

Write $\Z_{1/p}$ for the ring extension of $\Z$ by a single fraction $1/p$. 

\begin{cor}\label{cor:CohomIneq}
	Given a torus action $T\curvearrowright X$, fix a prime $p$ and assume $H^*(X,\Z)$ has no $q$-torsions for any prime $q\neq p$ (or no torsions at all), then the equality 
	$
	\sum_i \mathrm{rank}_{\Z}\,H^i(X^T,{\Z}) =\sum_i \mathrm{rank}_{\Z}\,H^i(X,{\Z})
	$
	holds if and only if the action $T\curvearrowright X$ is equivariantly formal in ${\Z_{1/p}}$ coefficients ($\Z$ coefficients). Moreover, if the equality does hold, then $H^*(X^T,\Z)$ has no $q$-torsions for any prime $q\neq p$ (or no torsions at all).
\end{cor}
\begin{proof}
	Suppose the action $T\curvearrowright X$ is equivariantly formal in ${\Z_{1/p}}$ coefficients. After tensoring the cohomology with $\Q$, we see the torus action is also equivariantly formal in $\Q$ coefficients and hence has the equality $\sum_i \dim_\Q\,H^i(X^T,\Q) = \sum_i \dim_\Q\,H^i(X,\Q)$ which gives us the equality regarding the $\Z$-ranks.
	
	Conversely, suppose we have the equality $\sum_i \mathrm{rank}_\Z\,H^i(X^T,\Z) = \sum_i \mathrm{rank}_\Z\,H^i(X,\Z)$. Consider the coefficients in $\K=\Q$ or $\F_q$ of any prime $q\neq p$. Since the free $\Z$-summands of $H^*(X^T,\Z)$ contribute to the $\K$-summands of $H^*(X^T,\K)$, we have the inequality $\sum_i \dim_\K\,H^i(X^T,\K) \geq \sum_i \mathrm{rank}_\Z\,H^i(X^T,\Z)$. Since we assume $H^*(X,\Z)$ has no $q$-torsions, all the $\K$-summands of $H^*(X,\K)$ exactly come from the free part of $H^*(X,\Z)$ and we have the equality $\sum_i \mathrm{rank}_\Z\,H^i(X,\Z) = \sum_i \dim_\K\,H^i(X,\K)$. Combining these three (in)equalities, we have $\sum_i \dim_\K\,H^i(X^T,\K) \geq \sum_i \dim_\K\,H^i(X,\K)$ which forces the inequality in Theorem~\ref{thm:CohomIneq} to be an equality, and implies the equivariant formality in coefficients $\Q$ or $\F_q$ of any prime $q\neq p$, equivalently in coefficient $\Z_{1/p}$. 
	
	As in the above argument, if the equality $\sum_i \mathrm{rank}_\Z\,H^i(X^T,\Z) = \sum_i \mathrm{rank}_\Z\,H^i(X,\Z)$ does hold, it forces the inequality $\sum_i \dim_\K\,H^i(X^T,\K) \geq \sum_i \mathrm{rank}_\Z\,H^i(X^T,\Z)$ to be an equality in coefficients $\K=\Q$ or $\F_q$ of any prime $q\neq p$. Hence, any $q$-torsion in $H^*(X^T,\Z)$ is ruled out.
\end{proof}

The fixed-point set can reveal information about the equivariant cohomology ring of the original manifold.
\begin{thm}[The Borel localization theorem, see \cite{AP93} p.\,133, Thm\,3.1.6]\label{thm:BorelLocal}
	Under the same assumptions as in Theorem\,\ref{thm:CohomIneq}, consider the inclusion $\iota:X^G\hookrightarrow X$, then the pullback homomorphism $\iota^*:H^*_G(X,\K)\rightarrow H^*_G(X^G,\K)$ is an isomorphism modulo $H^*_G(pt,\K)$-torsions. Moreover, if the action $G\curvearrowright X$ is equivariantly formal, then $\iota^*$ is injective.
\end{thm}

\subsubsection{Equivariant characteristic classes}

Given a rank-$n$ real vector bundle $\xi\rightarrow X$, if a group $G$ acts on both $\xi$ and $X$ and preserves the vector bundle structure, then $\xi\rightarrow X$ is a \textit{$G$-equivariant real vector bundle} with the \textit{$i$-th $G$-equivariant Stiefel-Whitney class} $w^G_i(\xi\rightarrow X)\triangleq w_i(EG\times_G\xi \rightarrow EG\times_G X)\in H^{i}_G(X,\F_2)$, the \textit{mod-$2$ $G$-equivariant Euler class} $e^G_{\F_2}(\xi\rightarrow X)\triangleq e_{\F_2}(EG\times_G\xi \rightarrow EG\times_G X)\in H^n_G(X,\F_2)$ and the \textit{$i$-th $G$-equivariant Pontryagin class} $p^G_i(\xi\rightarrow X)\triangleq p_i(EG\times_G\xi \rightarrow EG\times_G X)\in H^{4i}_G(X,\Z)$. If $\xi$ has an orientation which is preserved by the $G$-action, then we also have the \textit{integral $G$-equivariant Euler class} $e_{\Z}^G(\xi\rightarrow X)\triangleq e_{\Z}(EG\times_G\xi \rightarrow EG\times_G X)\in H^n_G(X,\Z)$.

In this paper, we will mostly consider the actions of tori or $2$-tori. Given a rank-$n$ torus $T^n$ or $2$-torus $\Z_2^n$, for the trivial actions of  $T^n$ or $\Z_2^n$ on a single point $pt$, we denote their equivariant cohomology by 
\begin{align*}
H^*_{T^n}(pt,\Z)&=H^*(BT^n,\Z)=H^*((\C P^\infty)^n,\Z)=\Z[\alpha_1,\ldots,\alpha_n]\\
H^*_{\Z_2^n}(pt,\F_2)&=H^*(B\Z_2^n,\F_2)=H^*((\R P^\infty)^n,\F_2)=\F_2[\beta_1,\ldots,\beta_n]
\end{align*}
where $\alpha_i$'s of cohomological degree $2$ are the negative generators of $H^*((\C P^\infty)^n,\Z)$, and $\beta_i$'s of cohomological degree $1$ are the generators of $H^*((\R P^\infty)^n,\F_2)$.  

\begin{lem}\label{lem:equivCharPt}
	Let $\tilde{\R}$ and $\tilde{\R}^2$ be the standard real representations of $\Z_2$ and $S^1$ by reflections and rotations respectively. They can also be viewed as $\Z_2$ and $S^1$-equivariant bundles over a single point. Write $H^*_{\Z_2}(pt,\F_2)=\F_2[\beta]$ and $H^*_{S^1}(pt,\Z)=\Z[\alpha]$. We have the equivariant characteristic classes $w^{\Z_2}(\tilde{\R})=1+\beta,\,e^{\Z_2}_{\F_2}(\tilde{\R})=\beta,\,
	p^{S^1}(\tilde{\R}^2)=1+\alpha^2,\,
	e_{\Z}^{S^1}(\tilde{\R}^2)=\alpha$.
\end{lem}
\begin{proof}
	The Borel construction $S^\infty \times_{\Z_2} \tilde{\R}$ is the universal real line bundle over the classifying space $\R P^\infty$. We have $w^{\Z_2}(\tilde{\R})=w(S^\infty \times_{\Z_2} \tilde{\R})=1+\beta$ and $e_{\F_2}^{\Z_2}(\tilde{\R})=w_1^{\Z_2}(\tilde{\R})=\beta$. Similarly, let $S^1\curvearrowright \tilde{\C}$ be the standard complex representation of rotations, the Borel construction $S^\infty \times_{S^1} \tilde{\C}$ is the universal complex line bundle over the classifying space $\C P^\infty$. Identifying $\tilde{\C} = \tilde{\R}^2$, we then have $e^{S^1}(\tilde{\R}^2)=e^{S^1}(\tilde{\C})=e(S^\infty \times_{S^1}\tilde{\C})=\alpha$. The square formula of Euler class gives $p^{S^1}_1(\tilde{\R}^2)=\big(e^{S^1}(\tilde{\R}^2)\big)^2=\alpha^2$. Hence $p^{S^1}(\tilde{\R}^2)=1+p^{S^1}_1(\tilde{\R}^2)=1+\alpha^2$.
\end{proof}

The canonical torus or $2$-torus actions on the real or oriented flag manifolds extend naturally to the universal vector bundles and make those bundles equivariant. 

Under the same notations of Subsubsection\,\ref{subsubsec:canTorusAction}, the $\Z_2^N$-action on $Fl(n_1,\ldots,n_k,\R^n)$ is induced from the $\Z_2^N$-action on $\R^n=\R^N\oplus \R^{n-N}$ where $\R^N$ is the standard representation of $\Z_2^N$ and $\R^{n-N}$ is a trivial representation. The $T^M$-action on $Fl(n_1,\ldots,n_k,\R^n)$ is induced from the $T^M$-action on $\R^n=\R^{2M}\oplus \R^{n-2M}$ where $\R^M$ is the standard representation of $T^M$ and $\R^{n-2M}$ is a trivial representation. Besides, simple verification shows that the topologically trivial line bundles $\mathcal{L}_i$ are also equivariantly trivial.

Denote the total equivariant Stiefel-Whitney, Pontryagin classes and the equivariant integral Euler class of the $j$-th universal bundle $\mathcal{W}_j$ by $\tilde{w}^{(j)}=1+\tilde{w}_1^{(j)}+\cdots+\tilde{w}_{n_j}^{(j)}$, $\tilde{p}^{(j)}=1+\tilde{p}_1^{(j)}+\cdots+\tilde{p}_{[n_j/2]}^{(j)}$ and $\tilde{e}$. Applying the Whitney product formulas, the square formula of Euler class and  Lemma~\ref{lem:equivCharPt} to decompositions\,\eqref{eq:fullDecom}\,\eqref{eq:completeDecom}, the equivariant versions of Relations\,\eqref{eq:Whitney}\,\eqref{eq:sqEuler}\,\eqref{eq:WhitneyEuler} are
\begin{gather}
\prod_{j=1}^{k}\tilde{w}^{(j)}=\prod_{l=1}^{N}(1+\beta_l), \qquad \qquad \qquad  \prod_{j=1}^{k}\tilde{p}^{(j)}=\prod_{l=1}^{M}(1+\alpha^2_l) \label{eq:EquivWhitney} \\
(\tilde{e}^{(j)})^2=\tilde{p}^{(j)}_{[n_j/2]}\label{eq:EquivsqEuler}\\
\prod_{j=1}^{k}\tilde{e}^{(j)}=\tilde{e}(\mathcal{V}) \label{eq:EquivWhitneyEuler}
\end{gather}
where $\tilde{e}(\mathcal{V})$, the equivariant Euler class of the full universal bundle $\mathcal{V}$, will be specified shortly in later sections.

\vskip 20pt
\section{The main tools and strategy}
\vskip 15pt
The cohomology rings of real Stiefel manifolds and real Grassmannians are well known in the coefficient rings $\F_2$ and $\Z_{1/2}$ which is the ring extension of $\Z$ by the fraction $1/2$. As introduced in Subsubsection\,\ref{subsubsec:generalFlag}, a real flag manifold is a bundle with the base a real Stiefel manifold and the fibre a usual real flag manifold. This bundle structure and its variations, together with the theorems of Leray and Borel will enable us to describe the cohomology rings of a real flag manifold in $\F_2$ and $\Z_{1/2}$ coefficients.  

\subsection{Mod-$2$ and $\Z_{1/2}$-coefficient ordinary cohomology rings of the real Stiefel manifolds}
For the real Stiefel manifold $V_{n-k}(\R^n)=SO(n)/SO(k)=O(n)/O(k)$, Ehresmann\,\cite{Eh39} introduced a cellular decomposition and announced the results of homology groups. Borel\,\cite{Bo53a} and Miller\,\cite{Mi53} described the cohomology rings in $\Z_{1/2}$ and $\F_2$ coefficients (also see the notes and books \cite{MT91,Fu04,BCM}).
\begin{thm}[\cite{Bo53a,Mi53}]\label{thm:realStief}
	All torsions in the integral cohomology of a real Stiefel manifold are of order $2$.
	\begin{enumerate}[(a)]
		\item For $\F_2$ coefficients, the cohomology ring of the real Stiefel manifold is 
		\[
		H^*\Big(\frac{O(n)}{O(k)},\F_2\Big)= \frac{\F_2[h_{i}\mid_{k\leq i \leq n-1}]}{\big\langle h_{i}^2=h_{2i} \mid_{k\leq i \leq n-1}\big\rangle}
		\]
		where $h_i$ is of degree $i$, and the denominator is the ideal generated by the relations $h_{i}^2=h_{2i}$ for $k\leq i \leq n-1$ and $h_{2i}$ is understood to be $0$ if $2i>n-1$. There is an additive basis of $H^*\big(O(n)/O(k),\F_2\big)$:
		\[
		h_{i_1}\ldots h_{i_l}, \qquad k\leq i_1<\cdots<i_l\leq n-1.
		\]
		\item For $\Z_{1/2}$ coefficients, the cohomology rings of the real Stiefel manifolds are exterior algebras
		\begin{align*}
		H^*(SO(2m)/SO(2k+1),\Z_{1/2})&=\Lambda_{\Z_{1/2}}[y_{k+1},\ldots,y_{m-1},x_m]\\
		H^*(SO(2m+1)/SO(2k+1),{\Z_{1/2}})&=\Lambda_{\Z_{1/2}}[y_{k+1},\ldots,y_{m}]\\
		H^*(SO(2m)/SO(2k),{\Z_{1/2}})&=\Lambda_{\Z_{1/2}}[e_{k},y_{k+1},\ldots,y_{m-1},x_m]\\
		H^*(SO(2m+1)/SO(2k),{\Z_{1/2}})&=\Lambda_{\Z_{1/2}}[e_{k},y_{k+1},\ldots,y_{m}]
		\end{align*}
		where $x_m,y_i$ are of odd degrees $2m-1,4i-1$, and $e_k$ of even degree $2k$ is the Euler class of the universal bundle $\mathcal{V}$ over $SO(2m)/SO(2k)$ or $SO(2m+1)/SO(2k)$.
	\end{enumerate} 
\end{thm}

\begin{rmk}
	If $k>0$, then $O(n)/O(k)=SO(n)/SO(k)$. We prefer using $O(n)/O(k)$ when discussing $2$-torus actions and mod-$2$ cohomology, while prefer $SO(n)/SO(k)$ when discussing torus actions and cohomology in other coefficients.
\end{rmk}

Consider a series of $2$-tori $\Z_2^l\subset O(l)\subset O(k)$ on the upper left $l$-by-$l$ block of $O(k)$ where $l\leq k$ and similarly the tori $T^{l'}\subset SO(2l')\subset SO(k)$ where $2l'\leq k$. Let them act on the real Stiefel manifold $O(n)/O(k)=SO(n)/SO(k)$ by left multiplications. Similar to Lemmas\,\ref{lem:fixedPointOfRealGeneralFlag},\,\ref{lem:fixedPointOfRealGeneralFlag2} on the fixed-point sets of real flag manifolds, we have

\begin{lem}\label{lem:fixedPointOfRealStiefel}
	The $\Z_2^l$-fixed-point set of $O(n)/O(k)$ is $O(n-l)/O(k-l)$, and the $T^{l'}$-fixed-point set of $SO(n)/SO(k)$ is $SO(n-2l')/SO(k-2l')$, where $O(n-l)$ and $SO(n-2l')$ are on the lower right blocks of $O(n)$.
\end{lem}

\begin{prop}\label{prop:FormalStiefel}
	The $\Z_2^l$-action on $O(n)/O(k)$, where $l\leq k$, is equivariantly formal in $\F_2$ coefficients and the $T^{l'}$-action on $SO(n)/SO(k)$, where $2l'\leq k$, is equivariantly formal in $\Z_{1/2}$ coefficients.
\end{prop}
\begin{proof}
	By Theorem\,\ref{thm:realStief}\,(a), the mod-$2$ Poincare polynomials of $O(n)/O(k)$ and its $\Z_2^l$-fixed-point set $O(n-l)/O(k-l)$ are $\prod_{i=k}^{n-1}(1+t^{i})$ and $\prod_{i=k-l}^{n-l-1}(1+t^{i})$. Evaluated at $t=1$, both polynomials give the same mod-$2$ total Betti number $2^{n-k}$. By Theorem\,\ref{thm:CohomIneq}, the $\Z_2^l$-action is equivariantly formal in $\F_2$ coefficients. Similarly, using Theorem\,\ref{thm:realStief}\,(b) and Corollary\,\ref{cor:CohomIneq}, we get the equivariant formality for the $T^{l'}$-actions in $\Z_{1/2}$ coefficients. 
\end{proof}

\subsection{The Leray-Borel theorems}
When a compact connected Lie group $G$ and its closed connected subgroup $H$ contain the same maximal torus, the rational cohomology ring of $G/H$ was described by Leray\,\cite{Le50}. If $H^*(G,\Z)$ and $H^*(H,\Z)$ have no $p$-torsions for a prime $p$ or no torsions at all, Borel\,\cite{Bo53a} showed that the description of the cohomology ring of $G/H$ also holds for $\F_p$ or $\Z$ coefficients. 

The original discussions by Leray and Borel have a more general form in terms of principal bundles and their quotient bundles. Given a principal $G$-bundle $G\rightarrow E \rightarrow E/G$, we automatically have a principal $H$-bundle $H\rightarrow E \rightarrow E/H$ over $E/H$, whose base manifold $E/H$ is also a bundle $G/H\rightarrow E/H \rightarrow E/G$ which we will call an \textit{associated bundle} of the principal $G$-bundle $G\rightarrow E \rightarrow E/G$.

Let $T$ be a maximal torus of $G$. On the quotient $E/T$, there is a well-defined right action of the normalizer $N_G(T)$, which descends to an action of the Weyl group $W_G$ and further induces a $W_G$-action on the cohomology ring $H^*(E/T,\Q)$.

\begin{thm}[The Leray-Borel theorem for associated bundles,\,\cite{Le50,Bo53a}]\label{thm:BundleLerayBorel}
	For a compact connected Lie group $G$ and its closed connected subgroup $H$ that contains the same maximal torus $T$, let $E$ be a principal $G$-bundle, then we have
	\begin{enumerate}[label=(\alph*)]
		\item the projection $E/T\rightarrow E/G$ induces an injective ring homomorphism $H^*(E/G,\Q)\hookrightarrow H^*(E/T,\Q)$ with the image $H^*(E/T,\Q)^{W_G}$, i.e. as a ring
		\[
		H^*(E/G,\Q) = H^*(E/T,\Q)^{W_G},
		\]
		\item the Leray-Serre spectral sequence of the bundle $G/H\rightarrow E/H \rightarrow E/G$ collapses at the $E_2$ page with an $H^*(E/G,\Q)$-module isomorphism
		\[
		H^*(E/H,\Q) \cong H^*(E/G,\Q)\otimes_\Q H^*(G/H,\Q)
		\]
		which can be modified to be an $H^*(E/G,\Q)$-algebra identification
		\[
		H^*(E/H,\Q)= H^*(E/G,\Q)\otimes_{H^*(BG,\Q)} H^*(BH,\Q),
		\]
		\item the Poincare polynomials of $E/H$, $G/H$ and $E/G$ satisfy the relation
		\[
		P(E/H,t)=P(E/G,t)\cdot P(G/H,t).
		\]
	\end{enumerate}
Moreover, if $H^*(G,\Z)$ and $H^*(H,\Z)$ have no $p$-torsions for a prime $p$ or no torsions at all, then the above description of the cohomology ring also holds for $\F_p$ or $\Z$ coefficients. 
\end{thm}

\begin{rmk}\label{rmk:BundleLerayBorel}
	We gather some well-known results as direct consequences of the above theorem.
	\begin{enumerate}
		\item Conclusion (a) is usually referred as the splitting principle. Let $E=EG$, then (a) implies $H^*(BG,\Q)= H^*(BT,\Q)^{W_G}$ which is a polynomial ring by a theorem of Chevalley. For instance, the classical groups $U(n),\,Sp(n),\,SO(2n+1),\,SO(2n)$ all contain a maximal torus $T^n$. As subrings of $H^*(BT^n,\Z)=\Z[t_1,\ldots,t_n]$ where $t_i$'s are the degree-$2$ negative generators, 
		\begin{align*}
		H^*(BU(n),\Q)&=\Q[c_1,\ldots,c_n] &
		H^*(BSp(n),\Q)&=\Q[q_1,\ldots,q_n]\\
		H^*(BSO(2n+1),\Q)&=\Q[p_1,\ldots,p_n] &
		H^*(BSO(2n),\Q)&=\Q[p_1,\ldots,p_{n-1},e_n]
		\end{align*}
		where the $i$-th Chern class $c_i$ is the $i$-th elementary symmetric function in $t_1,\ldots,t_n$; the $i$-th quaternionic and Pontryagin classes $q_i,p_i$ are the $i$-th elementary symmetric functions in $t^2_1,\ldots,t^2_n$; and the Euler class $e_n=t_1\cdots t_n$. For $BU(n),BSp(n)$, the coefficient ring can be improved to be $\Z$. for $BSO(2n+1),BSO(2n)$, the coefficient ring can be improved to be $\Z_{1/2}$.
		\item Let $E=G$ and assume $G,H$ are of the same rank, (b) implies $H^*(G/H,\Q)= \Q\otimes_{H^*(BG,\Q)} H^*(BH,\Q)$. The right hand side can also be written as $H^*(BH,\Q)/\langle H^+(BG,\Q)\rangle$
		where $\langle H^+(BG,\Q)\rangle$ is the ideal in $H^*(BH,\Q)= H^*(BT,\Q)^{W_H}$ generated by the positive-degree elements of $H^*(BG,\Q)= H^*(BT,\Q)^{W_G}$.
		\item If another compact Lie group $K$ acts on $E$ and preserves the principal $G$-bundle structure, then the Borel construction $EK\times_K E\rightarrow EK\times_K E/G$ is also a principal $G$-bundle. Substituting $EK\times_K E$ for the $E$ in the above Leray-Borel theorem, we then get a $K$-equivariant version of the Leray-Borel Theorem\,\ref{thm:BundleLerayBorel} after replacing the ordinary cohomology by the $K$-equivariant cohomology, and get the equivalence between the equivariant formality of $K\curvearrowright E/G$ and that of $K\curvearrowright E/T$. 
		\item Given an action $G\curvearrowright X$, we can substitute the product $EG\times X$ for the $G$-bundle $E$ in (a) and get $H_G^*(X,\Q)= H_T^*(X,\Q)^{W_G}$. Together substituting $EG\times X$ in (b), we see the equivariant formality of $G\curvearrowright X$ in rational coefficients is equivalent to that of $T\curvearrowright X$.
	\end{enumerate}
\end{rmk}

\begin{rmk}
	For a homogeneous space $G/H$ satisfying the equal-rank condition, various generalizations of the Leray-Borel-type description $H^*(G/H,\K)=H^*(BT,\K)^{W_H}/\langle H^+(BT,\K)^{W_G}\rangle$ to other ordinary or equivariant cohomology theories or for a Kac-Moody group $G$ have appeared in \cite{KK86,KK90,Ar89,KV93,Br97,KiKr13,CZZ15}.
\end{rmk}

\begin{rmk}
	Baum, Smith\,\cite{BS67} and Wolf\,\cite{Wo78} studied the cohomology of the bundle $E/H$ using the Eilenberg-Moore spectral sequences. They didn't require $H$ to be of the same rank as $G$ but assumed the base $E/G$ to be of certain types of homogeneous spaces. 
\end{rmk}

\begin{rmk}
	Jeff Carlson kindly pointed out the work of Leray to the author. There is an extensive discussion on cohomology of homogeneous spaces in his thesis\,\cite{Carl15}.
\end{rmk}

Similarly, Borel considered the maximal $2$-torus $\Z_2^n$ of a Lie group $G$ and the mod-$2$ cohomology.
\begin{thm}[The mod-$2$ Leray-Borel theorem for associated bundles,\,\cite{Bo53b}]\label{thm:BundleLerayBorel2}
	For a compact, possibly disconnected Lie group $G$ and its closed subgroup $H$ containing the same maximal $2$-torus $\Z_2^n$, suppose the inclusions $\Z_2^n\subseteq H\subseteq G$ induce injective pullbacks $H^*(BG,\F_2)\hookrightarrow H^*(BH,\F_2)\hookrightarrow H^*(B\Z_2^n,\F_2)$ and there are ring isomorphisms
	\[
	H^*(G/\Z^n_2,\F_2)\cong\F_2\otimes_{H^*(BG,\F_2)} H^*(B\Z_2^n,\F_2) \qquad \qquad H^*(H/\Z^n_2,\F_2)\cong\F_2\otimes_{H^*(BH,\F_2)} H^*(B\Z_2^n,\F_2).
	\]
	Let $E$ be a compact $G$-bundle, then we have
	\begin{enumerate}[label=(\alph*)]
		\item the Leray-Serre spectral sequence of the bundle $G/H\rightarrow E/H \rightarrow E/G$ collapses at the $E_2$ page with an $H^*(E/G,\F_2)$-module isomorphism
		\[
		H^*(E/H,\F_2) \cong H^*(E/G,\F_2)\otimes_{\F_2} H^*(G/H,\F_2)
		\]
		which can be modified to be an $H^*(E/G,\F_2)$-algebra identification
		\[
		H^*(E/H,\F_2)= H^*(E/G,\F_2)\otimes_{H^*(BG,\F_2)} H^*(BH,\F_2),
		\]
		\item the mod-$2$ Poincare polynomials of $E/H$, $G/H$ and $E/G$ satisfy the relation
		\[
		P_2(E/H,t)=P_2(E/G,t)\cdot P_2(G/H,t).
		\]
	\end{enumerate}
\end{thm}

Similar to Remark\,\ref{rmk:BundleLerayBorel}, we have
\begin{rmk}\label{rmk:BundleLerayBorel2}
	Some direct consequences of the above theorem are:
	\begin{enumerate}
		\item The classical groups $O(n),\,T^n=U(1)^n,\,U(n),\,Sp(n)$ all contain the $2$-torus $O(1)^n=\Z_2^n$. Borel showed that, as subrings of $H^*(B\Z_2^n,\F_2)=\F_2[s_1,\ldots,s_n]$ where $s_i$'s are the degree-$1$ generators,
		\begin{align*}
		H^*(BO(n),\F_2)&=\F_2[w_1,\ldots,w_n] &
		H^*(BT^n,\F_2)&=\F_2[t_1,\ldots,t_n]\\
		H^*(BU(n),\F_2)&=\F_2[c_1,\ldots,c_n] &
		H^*(BSp(n),\F_2)&=\F_2[q_1,\ldots,q_n]
		\end{align*}
		where the $i$-th Stiefel-Whitney, Chern and quaternionic classes $w_i,c_i,q_i$ are the $i$-th elementary symmetric functions in $s_1,\ldots,s_n;s^2_1,\ldots,s^2_n;s^4_1,\ldots,s^4_n$ respectively, and $t_i=s_i^2$.
		\item Let $E=G$, then (b) implies $H^*(G/H,\F_2)= \F_2\otimes_{H^*(BG,\F_2)} H^*(BH,\F_2)= H^*(BH,\F_2)/\langle H^+(BG,\F_2)\rangle$. 
		\item If another $2$-torus $\Z_2^l$ acts on $E$ and preserves the principal $G$-bundle structure, then we get a $\Z_2^l$-equivariant version of the mod-$2$ Leray-Borel Theorem\,\ref{thm:BundleLerayBorel2} after replacing the ordinary cohomology by the $\Z_2^l$-equivariant cohomology, and get the equivalence between the equivariant formality of $\Z_2^l\curvearrowright E/G$ and that of $\Z_2^l\curvearrowright E/T$. 
	\end{enumerate}
\end{rmk}

\subsection{The main strategy}
Let $G$ be a compact connected Lie group, and $H$ be a closed, possibly disconnected subgroup. Our goal is to understand the cohomology ring of the homogeneous space $G/H$.
\begin{description}
	\item[\textit{Step 1}] If $H$ is disconnected, let $H_0$ be its identity component, we can consider the covering
	$
	H/H_0 \rightarrow G/H_0 \rightarrow G/H
	$
	and use a covering-space argument to get
	$
	H^*(G/H,\K)= H^*(G/H_0,\K)^{H/H_0}
	$.
	Then we switch to $G/H_0$.
	\item[\textit{Step 2}] Let $T$ be a maximal torus of $H_0$, we consider the bundle
	$
	H_0/T \rightarrow G/T \rightarrow G/H_0
	$
	and use the Leray-Borel Theorem\,\ref{thm:BundleLerayBorel}\,(a) to get
	$
	H^*(G/H_0,\K)\cong H^*(G/T,\K)^{W_{H_0}}.
	$
	Then we switch to $G/T$.
	\item[\textit{Step 3}] If $T$ is also a maximal torus of $G$, then use Remark\,\ref{rmk:BundleLerayBorel} to get
	$
	H^*(G/T,\K)= H^*(BT,\K)/\langle H^+(BG,\K)\rangle.
	$
	\item[\textit{Step 4}]
	Otherwise, $G/T$ might still be not easy to understand, then we try a reverse of the above step. Let $K$ be a closed connected subgroup of $G$ and contain $T$ as a maximal torus, we consider the bundle
	$
	K/T \rightarrow G/T \rightarrow G/K
	$
	and use the Leray-Borel Theorem\,\ref{thm:BundleLerayBorel} to get
	$
	H^*(G/T,\K) = H^*(G/K,\K)\otimes_{H^*(BK,\K)} H^*(BT,\K)
	$
	where $H^*(BK,\K)= H^*(BT,\K)^{W_K}$. Then we switch to $G/K$, which hopefully can be understood using some other methods (Cartan models, Leray-Serre spectral sequences, Eilenberg-Moore spectral sequences, etc.).
\end{description}

\begin{framed}
	\bf Roughly speaking, if we already know the cohomology of $G/K$, then we will know the cohomology of $G/H$ for any $H$ that shares the same maximal torus of $gKg^{-1}$ for some $g\in G$.
\end{framed}

In applications, we could repeat some steps several times. The coefficient ring $\K$ is usually $\Q$ but could be improved if the integral cohomology of $G,H,K$ does not have certain torsions. Using the mod-$2$ Leray Borel theorem\,\ref{thm:BundleLerayBorel2} and considering the maximal $2$-torus instead, the above strategy also works for $\F_2$ coefficients in many cases. 

\vskip 20pt
\section{Mod-$2$ equivariant cohomology rings of the real flag manifolds}
\vskip 15pt
When $n=n_1+\cdots+n_k$, both $O(n)$ and $O(n_1)\times \cdots\times O(n_k)$ contain the same maximal $2$-torus $O(1)^n=\Z_2^n$. Borel showed that
\begin{thm}[\cite{Bo53b}]
	If $n=n_1+\cdots+ n_k$, as a ring,
	\begin{align*}
	H^*\Big(\frac{O(n)}{O(n_1)\times \cdots\times O(n_k)},\F_2\Big)=\F_2\otimes_{H^*(BO(n),\F_2)} H^*(BO(n_1)\times \cdots\times BO(n_k),\F_2)=\frac{\F_2[w_i^{(j)}\mid_{1\leq j\leq k,1\leq i\leq n_j}]}{\big\langle\prod_{j=1}^{k}w^{(j)}=1\big\rangle}
	\end{align*}
	where $w^{(j)}=1+w_1^{(j)}+\cdots+w_{n_j}^{(j)}$ is the total Stiefel-Whitney class of the $j$-th universal bundle over the real flag manifold. The mod-$2$ Poincare polynomial is
	\begin{align*}
	P_{2}\Big(\frac{O(n)}{O(n_1)\times \cdots\times O(n_k)},t\Big)=\frac{P_{2}(BO(n_1),t)\cdots P_{2}(BO(n_k),t)}{P_{2}(BO(n),t)}=\frac{\prod_{i=1}^{n} (1-t^{i})}{\prod_{j=1}^k\prod_{i=1}^{n_j}(1-t^{i})}.
	\end{align*}
\end{thm}

\begin{rmk}
	When $k>1$, the $2$-rank of $SO(n_1)\times \cdots\times SO(n_k)$ is less than that of $SO(n)$, hence the mod-$2$ Leray-Borel Theorem does not apply to the oriented flag manifolds.
\end{rmk}

The above expression of the Poincare polynomial can be written more succinctly using combinatorial notations.
\begin{dfn}
	Let $n_1,\ldots,n_k$ be positive integers and $n=n_1+\cdots+n_k$, define the \textit{$t$-number}, \textit{$t$-factorial} and \textit{Gaussian multinomial coefficient} progressively as 
	\begin{align*}
	[n]_{t} &\triangleq 1+t+\cdots+t^{n-1}=\frac{1-t^n}{1-t}\\
	[n]_{t}! &\triangleq [1]_{t}\cdots [n]_{t}\\
	\binom{n}{n_1,\ldots,n_k}_t &\triangleq \frac{[n]_{t}!}{[n_1]_{t}!\cdots[n_k]_{t}!}.
	\end{align*} 
\end{dfn}

The mod-$2$ Poincare polynomial of the real flag manifold then can be written as
\[
P_2\Big(\frac{O(n)}{O(n_1)\times \cdots\times O(n_k)},t\Big)=\binom{n}{n_1,\ldots,n_k}_t.
\]

\begin{exm}
	Consider the real Grassmannian of $k$-dimensional subspaces in $\R^n$ and the real complete flag manifold of flags $0\subset V_1 \subset V_2 \cdots \subset V_{n-1}  \subset \R^{n}$:
	{
		\begin{align*}
		G_k(\R^n)=Fl(k,n-k,\R^n)=\frac{O(n)}{O(k)\times O(n-k)} \qquad \qquad 
		Fl(\R^n)=Fl(\underbrace{1,\ldots,1}_{n \text{ items}},\R^n)=\frac{O(n)}{O(1)^n}.
		\end{align*}}The mod-$2$ cohomology rings (given by Chern\,\cite{Ch48} for $G_k(\R^n)$, and by Borel\,\cite{Bo53b} for $Fl(\R^n)$) are:
	{
		\[
		H^*(G_k(\R^n),\F_2)=\frac{\F_2[w_1,\ldots,w_k;w'_1,\ldots,w'_{n-k}]}{\big\langle ww'=1\big\rangle} \qquad \qquad H^*(Fl(\R^n),\F_2)=\frac{\F_2[s_1,\ldots,s_n]}{\big\langle\prod_{i=1}^{n}(1+s_i)=1\big\rangle}
		\]}where $w=1+w_1+\cdots w_k,w'=1+w'_1+\cdots w'_{n-k}$ are the total Stiefel-Whitney classes of the universal bundle and its complementary bundle over $G_k(\R^n)$, and the $s_i$ of degree $1$ is the first Stiefel-Whitney class of the $i$-th universal real line bundle over $Fl(\R^n)$. The Poincare polynomials are respectively
	{
		\begin{align*}
		\frac{\prod_{i=1}^{n} (1-t^{i})}{\prod_{i=1}^{k} (1-t^{i})\prod_{i=1}^{n-k} (1-t^{i})}=\binom{n}{k,n-k}_{t}, \qquad \qquad
		\frac{\prod_{i=1}^{n} (1-t^{i})}{(1-t)^n}=\binom{n}{1,\ldots,1}_{t}=[n]_{t}!.
		\end{align*}}
\end{exm}

\subsection{Reduction from real flag manifolds to real Stiefel manifolds}
Only assuming $n\geq N= n_1+\cdots+ n_k$, we shall consider the bundle 
\[
\frac{O(N)}{O(n_1)\times \cdots \times O(n_k)}
\longhookrightarrow
\frac{O(n)}{O(n_1)\times \cdots \times O(n_k)}
\longrightarrow
\frac{O(n)}{O(N)}
\]
which is equivariant under the canonical actions of the $2$-torus $O(1)^N=\Z_2^N$.

By Remark\,\ref{rmk:BundleLerayBorel2}\,(3), we can apply a $\Z_2^N$-equivariant version of the mod-$2$ Leray-Borel Theorem\,\ref{thm:BundleLerayBorel2} to the above equivariant bundle, and then get $H^*_{\Z_2^N}({O(n)}/{O(N)},\F_2)$-module isomorphisms:
\begin{align*}
H^*_{\Z_2^N}\Big(\frac{O(n)}{O(n_1)\times\cdots \times O(n_k)},\F_2\Big)
&\cong H^*_{\Z_2^N}\Big(\frac{O(n)}{O(N)},\F_2\Big)\otimes_{\F_2}
H^*\Big(\frac{O(n)}{O(n_1)\times \cdots\times O(n_k)},\F_2\Big)\\
&\cong H^*_{\Z_2^N}\Big(\frac{O(n)}{O(N)},\F_2\Big)\otimes_{\F_2}  \frac{\F_2[w_i^{(j)}\mid_{1\leq j\leq k,1\leq i\leq n_j}]}{\big\langle\prod_{j=1}^{k}w^{(j)}=1\big\rangle}.
\end{align*}

We can further lift the ordinary Stiefel-Whitney classes $w_i^{(j)}$ to be the $\Z_2^N$-equivariant Stiefel-Whitney classes $\tilde{w}_i^{(j)}$ and lift the relation $\prod_{j=1}^{k}w^{(j)}=1$ in Relation\,\eqref{eq:Whitney} to its $\Z_2^N$-equivariant version in Relation\,\eqref{eq:EquivWhitney}, then get:

\begin{prop}\label{prop:mod2EquivCohRealFlag}
	As an $H^*_{\Z_2^N}\big({O(n)}/{O(N)},\F_2\big)$-algebra, 
	\begin{align*}
	H^*_{\Z_2^N}\Big(\frac{O(n)}{O(n_1)\times\cdots \times O(n_k)},\F_2\Big)
	=  \frac{H^*_{\Z_2^N}\big(\frac{O(n)}{O(N)},\F_2\big)[\tilde{w}_i^{(j)}\mid_{1\leq j\leq k,1\leq i\leq n_j}]}{\big\langle\prod_{j=1}^{k}\tilde{w}^{(j)}=\prod_{l=1}^{N}(1+\beta_l)\big\rangle}
	\end{align*}
	where $\tilde{w}^{(j)}=1+\tilde{w}_1^{(j)}+\cdots+\tilde{w}_{n_j}^{(j)}$ is the $\Z_2^N$-equivariant total Stiefel-Whitney class of the $j$-th universal bundle over the real flag manifold. The $\Z_2^N$-action on ${O(n)}/{O(n_1)\times\cdots \times O(n_k)}$ is equivariantly formal, hence the mod-$2$ ordinary cohomology ring can be obtained by substituting $\beta_1=\cdots=\beta_N=0$ and the mod-$2$ ordinary Poincare polynomial is 
	\begin{align*}
	\prod_{i=N}^{n-1}(1+t^{i})\cdot \binom{N}{n_1,\ldots,n_k}_{t}.
	\end{align*}
\end{prop}
It remains to understand the ring structure of $H^*_{\Z_2^N}\big({O(n)}/{O(N)},\F_2\big)$ in the next subsection.

\subsection{The case of real Stiefel manifolds}By Proposition\,\ref{prop:FormalStiefel}, the $\Z_2^l$-actions on $O(n)/O(k)$ is equivariantly formal in $\F_2$ coefficients, where $l\leq k$. Hence we have $H_{\Z_2^l}^*(pt,\F_2)=\F_2[\beta_1,\ldots,\beta_l]$-module isomorphisms
\begin{align*}
H_{\Z_2^l}^*\Big(\frac{O(n)}{O(k)},\F_2\Big)\cong H_{\Z_2^l}^*(pt,\F_2)\otimes_{\F_2} H^*\Big(\frac{O(n)}{O(k)},\F_2\Big)
\cong \F_2[\beta_1,\ldots,\beta_l]\otimes_{\F_2}\frac{\F_2[h_{i}\mid_{k\leq i \leq n-1}]}{\big\langle h_{i}^2=h_{2i} \mid_{k\leq i \leq n-1}\big\rangle}.
\end{align*}

However, there are two ambiguities about these isomorphisms:
\begin{enumerate}
	\item The $\F_2[\beta_1,\ldots,\beta_l]$-module isomorphisms do not give the ring structure of $H_{\Z_2^l}^*\big({O(n)}/{O(k)},\F_2\big)$.
	\item The isomorphisms are not canonical, i.e. the equivariant lifts of the generators $h_{i}$ are not uniquely determined.
\end{enumerate}

Note the $\Z_2^l$-action is free on $O(n-l)/O(k-l)$. As an $\F_2[\beta_1,\ldots,\beta_l]$-algebra,
$$H_{\Z_2^l}^*\Big(\frac{O(n-l)}{O(k-l)},\F_2\Big)= H_{\Z_2^l}^*(pt,\F_2)\otimes_{\F_2} H^*\Big(\frac{O(n-l)}{O(k-l)},\F_2\Big)
= \F_2[\beta_1,\ldots,\beta_l]\otimes_{\F_2}\frac{\Z_2[h_{i}\mid_{k-l\leq i \leq n-l-1}]}{\big\langle h_{i}^2=h_{2i} \mid_{k-l\leq i \leq n-l-1}\big\rangle}.$$ 

Consider the inclusion of the $\Z_2^l$-fixed-point set: $\iota_{n,k,l}: O(n-l)/O(k-l) \hookrightarrow O(n)/O(k)$ and its pullback homomorphism of equivariant cohomology rings
$$\iota^*_{n,k,l}: H_{\Z_2^l}^*\Big(\frac{O(n)}{O(k)},\F_2\Big) \longrightarrow H_{\Z_2^l}^*\Big(\frac{O(n-l)}{O(k-l)},\F_2\Big).$$

By Theorem\,\ref{thm:BorelLocal}, the $\F_2[\beta_1,\ldots,\beta_l]$-algebra homomorphism $\iota^*_{n,k,l}$ is injective. Therefore, understanding the ring structure of $H_{\Z_2^l}^*\big({O(n)}/{O(k)},\F_2\big)$ is equivalent to identifying the image of $\iota^*_{n,k,l}$. Especially, we need to specify the $\iota^*_{n,k,l}$-images of the equivariant lifts of $h_i$'s in order to resolve the two ambiguities. 

Given a graded vector space $V^*$, let 
$$\Psi_k:V^*\longrightarrow V^{*\geq k}$$
be the projection that only preserves the components of degrees no less than $k$. 

\begin{thm}\label{thm:mod2EquivRealStiefel}
	$H_{\Z_2^l}^*\big({O(n)}/{O(k)},\F_2\big)$ can be identified, via the injective homomorphism $\iota^*_{n,k,l}$, as an $\F_2[\beta_1,\ldots,\beta_l]$-subalgebra of 
	\[
	\F_2[\beta_1,\ldots,\beta_l]\otimes_{\F_2}\frac{\F_2[h_{i}\mid_{k-l\leq i \leq n-l-1}]}{\big\langle h_{i}^2=h_{2i} \mid_{k-l\leq i \leq n-l-1}\big\rangle}
	\]
	generated on $\tilde{h}_i$'s for $k\leq i\leq n-1$ defined by the formula:
	\begin{align*}
	\iota^*_{n,k,l}\Big(\sum_{i=k}^{n-1} \tilde{h}_i\Big)=\Psi_k\Big(\prod_{j=1}^{l}(1+\beta_j)\cdot\sum_{i=k-l}^{n-l-1}h_i\Big).
	\end{align*}
\end{thm}

For later use of induction, we first prove the case when $n-k=1$, i.e. we have $O(n)/O(k)=S^k=\big\{(x_1,\ldots,x_{k+1})\in\R^{k+1}\mid \sum_{i=1}^{k+1}x_i^2=1\big\}$ acted by $\Z_2^l=\big\{(\epsilon_1,\ldots,\epsilon_l)\mid \epsilon_i = \pm 1\big\}$ where $l\leq k$ in the way that $(\epsilon_1,\ldots,\epsilon_l)\cdot (x_1,\ldots,x_{k+1})=(\epsilon_1 x_1,\ldots,\epsilon_l x_l,x_{l+1}\ldots,x_{k+1})$ flips $l$ coordinates. The $\Z_2^l$-fixed-point set is $O(k-l+1)/O(k-l)=S^{k-l}=\big\{(x_{l+1},\ldots,x_{k+1})\in 0^l\oplus \R^{k-l+1}\mid \sum_{i=l+1}^{k+1}x_i^2=1\big\}$. If $k-l\geq 1$, then $H^*(S^{k-l},\F_2)=\F_2[h_{k-l}]/\langle h_{k-l}^2 \rangle$. If $k-l=0$, then $H^*(S^{k-l},\F_2)=\F_2\oplus \F_2$. Write $h_0=(0,1)$, then $h_0^2=h_0$ and $H^*(S^0,\F_2)=\F_2[h_0]/\langle h_0^2-h_0\rangle$.

\begin{prop}\label{prop:mod2LocalizeSphere}
	Let $\iota:S^{k-l}\hookrightarrow S^{k}$ be the inclusion of fixed-point set. The image of the pullback monomorphism 
	\[
	\iota^*:H^*_{\Z^l_2}(S^{k},\F_2)\longhookrightarrow H^*_{\Z^l_2}(S^{k-l},\F_2)=
	\begin{cases}
	\F_2[\beta_1,\ldots,\beta_l][h_{k-1}]/\big\langle h_{k-l}^2\big\rangle, & k-l\geq 1\\
	\F_2[\beta_1,\ldots,\beta_l][h_0]/\big\langle h_0^2-h_0\big\rangle, & k-l=0,
	\end{cases}
	\]
    as an $\F_2[\beta_1,\ldots,\beta_l]$-subalgebra, is generated on $(\beta_1\cdots\beta_l) \cdot h_{k-l}$ which is exactly $\Psi_k\big((1+\beta_1)\cdots(1+\beta_l)h_{k-l}\big)$. 
\end{prop}
\begin{proof}
	The starting case is when $l=1$. $\Z_2\curvearrowright S^k$ flips the first coordinate with the fixed-point set $S^{k-1}$ where $k\geq 1$. Because of equivariant formality, the generator $h_{k}$ of $H^*(S^{k},\F_2)$ has an equivariant lift $\tilde{h}_{k}$ in $H^*_{\Z_2}(S^{k},\F_2)$. The image $\iota^*(\tilde{h}_{k})\in \F_2[\beta_1][h_{k-1}]/\langle h_{k-1}^2\rangle$ or $\F_2[\beta_1][h_{0}]/\langle h_{0}^2-h_{0}\rangle$, as a degree-$k$ element, is a linear combination of $\beta_1^{k},\beta_1 h_{k-1}$. Since $\tilde{h}_{k}$ is a nontrivial $\F_2[\beta_1]$-algebra generator, $\iota^*(\tilde{h}_{k})\neq \beta_1^{k}$. If $\iota^*(\tilde{h}_{k})=\beta_1^{k}+\beta_1 h_{k-1}$, we can replace $\tilde{h}_{k}$ by $\tilde{h}_{k}-\beta_1^{k}$ and manage to make $\iota^*(\tilde{h}_{k})=\beta_1 h_{k-1}$. 
	
	When $l\geq 2$, for any $1\leq i\leq l$, consider the $(k-l+1)$-spheres $S^{k-l+1}_i=\big\{(x_1,\ldots,x_{k+1})\in S^k \mid x_j=0 \mbox{ if } j\leq l \mbox{ and } j\neq i \big\}$ and the rank-$1$ $2$-sub-tori $\Z_{2,i}=\big\{(\epsilon_1,\ldots,\epsilon_l)\in \Z_2^l \mid \epsilon_j=1 \mbox{ if } j\neq i\big\}$. The sphere $S^{k-l+1}_i$ inherits a $\Z_2^l$-action from $S^k$. The rank-$1$ $2$-subtorus $\Z_{2,i}$ acts on $S^{k-l+1}_i$ by flipping the $x_i$-coordinate with the fixed-point set $S^{k-l}=\big\{(x_1,\ldots,x_{k+1})\in S^k \mid x_j=0 \mbox{ if } j\leq l\big\}$. We have proved the case of $l=1$ that the image of the monomorphism $\bar{\iota^*_i}: H^*_{\Z_{2,i}}(S^{k-l+1}_i,\F_2)\hookrightarrow H^*_{\Z_{2,i}}(S^{k-l}_i,\F_2)$, as an $H^*_{\Z_{2,i}}(pt,\F_2)=\F_2[\beta_i]$-subalgebra, is generated on $\beta_i h_{k-l}$. Since the other rank-$1$ $2$-subtori $\Z_{2,j}, j\neq i$ acts trivially, we have the ring isomorphisms $H^*_{\Z_{2}^l}(S^{k-l+1}_i,\F_2)=H^*_{\Z_{2}^l/\Z_{2,i}}(pt,\F_2)\otimes_{\F_2}H^*_{\Z_{2,i}}(S^{k-l+1}_i,\F_2)=\F_2[\beta_1,\ldots,\beta_{i-1},\beta_{i+1},\ldots,\beta_l]\otimes_{\F_2}H^*_{\Z_{2,i}}(S^{k-l+1}_i,\F_2)$. Hence, the $\F_2[\beta_i]$-generator $\beta_i h_{k-l}$ is also an $\F_2[\beta_1,\ldots,\beta_l]$-generator for the image of the monomorphism $\iota^*_i: H^*_{\Z_{2}^l}(S^{k-l+1}_i,\F_2)\hookrightarrow H^*_{\Z_{2}^l}(S^{k-l}_i,\F_2)$.
	
	Let $\tilde{h}_k\in H^*_{\Z^l_2}(S^{k},\F_2)$ be an equivariant lift of the degree-$k$ generator of $H^*(S^{k},\F_2)$. Via the monomorphism $\iota^*:H^*_{\Z^l_2}(S^{k},\F_2)\hookrightarrow H^*_{\Z^l_2}(S^{k-l},\F_2)=\F_2[\beta_1,\ldots,\beta_l]\otimes_{\F_2}H^*(S^{k-l},\F_2)$, we have $\iota^*(\tilde{h}_k)=f+g\cdot h_{k-l}$, where $f,g\in \F_2[\beta_1,\ldots,\beta_l]$ are homogeneous polynomials of degrees $k,l$. We can replace $\tilde{h}_k$ by $\tilde{h}_k-f$ to make $\iota^*(\tilde{h}_k)=g\cdot h_{k-l}$. Since $\tilde{h}_k$ is a nontrivial $\F_2[\beta_1,\ldots,\beta_l]$-generator, $g\neq 0$. The inclusion $\iota:S^{k-l}{\hookrightarrow} S^{k}$ factors through the inclusions $S^{k-l}\overset{\iota_i}{\hookrightarrow} S_i^{k-l+1}\hookrightarrow S^{k}$, so do their pullback homomorphisms of $\Z^l_2$-equivariant cohomology rings. We have $\mathrm{Im}\, \iota^*\subset \mathrm{Im}\, \iota^*_i \subset H^*_{\Z^l_2}(S^{k-l},\F_2)$. Since $\mathrm{Im}\, \iota^*_i$ is generated on $\beta_i h_{k-l}$, then $\iota^*(\tilde{h}_k)=g\cdot h_{k-l}$ is a multiple of $\beta_i h_{k-l}$ for $1\leq i\leq l$, hence a multiple of $(\beta_1\cdots\beta_l)h_{k-l}$, therefore $\iota^*(\tilde{h}_k)=(\beta_1\cdots\beta_l)h_{k-l}$ by degree reason. 
\end{proof}

\begin{rmk}
	The case of $\Z_2$ acting on $S^1$ by reflection has appeared in the work of Schmid\,\cite{Sc01}, and Biss, Guillemin and Holm\,\cite{BGH04}.
\end{rmk}

\begin{cor}
	For the $2$-torus actions on spheres, we have the equivariant cohomology rings:
	\begin{enumerate}[label=(\alph*)]
		\item If $\Z^l_2$ acts on $S^k$ by flipping $l$ coordinates, then
		\[
		H^*_{\Z^l_2}(S^{k},\F_2)=
		\begin{cases}
		\F_2[\beta_1,\ldots,\beta_l][\tilde{h}_{k}]/\big\langle \tilde{h}_{k}^2\big\rangle, & l< k\\
		\F_2[\beta_1,\ldots,\beta_k][\tilde{h}_{k}]/\big\langle \tilde{h}_{k}^2-(\beta_1\cdots\beta_k)\tilde{h}_{k}\big\rangle, & l= k\\
		\F_2[\beta_1,\ldots,\beta_k][\tilde{w}_1;\tilde{w}'_1,\ldots,\tilde{w}'_k]/\big\langle \tilde{w}\tilde{w}'=(1+\beta_1)\cdots(1+\beta_k)\big\rangle, & l= k+1.\\
		\end{cases}
		\]
		\item If $\Z_2$ acts on $S^k$ by flipping $l$ coordinates simultaneously, then
		\[
		H^*_{\Z_2}(S^{k},\F_2)=
		\begin{cases}
		\F_2[\beta][\tilde{h}_{k}]/\big\langle \tilde{h}_{k}^2\big\rangle, & l< k\\
		\F_2[\beta][\tilde{h}_{k}]/\big\langle \tilde{h}_{k}^2-\beta^k\tilde{h}_{k}\big\rangle, & l= k\\
		\F_2[\beta]/\big\langle \beta^{k+1}\big\rangle, & l= k+1.\\
		\end{cases}
		\]
	\end{enumerate}
\end{cor}

\begin{proof}
	Consider the rank-1 $2$-subtorus $H \triangleq \big\{(\epsilon_1,\ldots,\epsilon_{l}) \mid \epsilon_1=\cdots=\epsilon_{l} \in \{\pm 1\}\big\}\subset \Z_2^l$. The action of $H$ on $S^k$ is exactly the $\Z_2$-action that flips $l$ coordinates simultaneously.
	
	The first two cases of $\Z^l_2$-actions follow from Proposition\,\ref{prop:mod2LocalizeSphere} directly. The inclusion $H\subset \Z_2^l$ gives the pullback homomorphism $H^*_{\Z_2^l}(pt,\F_2)=\F_2[\beta_1,\ldots,\beta_l]\rightarrow H^*_{\Z_2}(pt,\F_2)=\F_2[\beta]$ by sending $\beta_i \mapsto \beta,\,\forall i$. Hence substituting $\beta_i=\beta,\,\forall i$ in the first two cases of $\Z^l_2$-action gives the first two cases of the simultaneous $\Z_2$-actions. 
	
	In the third cases when $l=k+1$, note that $H=\big\{(\epsilon_1,\ldots,\epsilon_{k+1}) \mid \epsilon_1=\cdots=\epsilon_{k+1} \in \{\pm 1\}\big\}$ acts on $S^k$ freely with the quotient $S^k/H=\R P^k$ on which there is an induced action of the quotient $2$-torus $\Z_2^{k+1}/H \cong \Z_2^{k}$. Then we have $H^*_{\Z_2^{k+1}}(S^k,\F_2)=H^*_{\Z_2^{k+1}/H}(S^k/H,\F_2)=H^*_{\Z_2^{k}}(\R P^k,\F_2)=H^*_{\Z_2^{k}}\big(O(k+1)/\big(O(k)\times O(1)\big),\F_2\big)$ which can be described by Proposition\,\ref{prop:mod2EquivCohRealFlag}. For the simultaneous $\Z_2$-action, $H^*_{\Z_2}(S^k,\F_2)=H^*(S^k/H,\F_2)=H^*(\R P^k,\F_2)$.
\end{proof}

For the general case when $1\leq l\leq k \leq n$, consider the following commutative diagram of Stiefel manifolds
\begin{center}
	\begin{tikzpicture}[description/.style={fill=white,inner sep=2pt}]
	\matrix (m) [matrix of math nodes, row sep=2em,
	column sep=2.5em, text height=1.5ex, text depth=0.25ex]
	{S^{k}=O(\R^{k+1})/O(\R^{k})  & O(\R^n)/O(\R^{k}) & O(\R^n)/O(\R^{k+1}) \\
	S^{k-l}=O(0^l \oplus\R^{k-l+1})/O(0^l \oplus\R^{k-l})  & O(0^l \oplus\R^{n-l})/O(0^l \oplus\R^{k-l}) & O(0^l \oplus\R^{n-l})/O(0^l \oplus\R^{k-l+1})\\};
	\path[->]
	(m-1-1) edge node[auto]{$j$} (m-1-2)
	(m-1-2) edge node[auto]{$\pi$} (m-1-3)
	(m-2-1) edge node[auto]{$j'$}(m-2-2)
	(m-2-2) edge node[auto]{$\pi'$} (m-2-3);
	\path[left hook->]
	(m-2-1) edge node[auto]{$\iota_{k+1,k,l}$} (m-1-1)
	(m-2-2) edge node[auto]{$\iota_{n,k,l}$} (m-1-2)
	(m-2-3) edge node[auto]{$\iota_{n,k+1,l}$} (m-1-3);
	\end{tikzpicture}
\end{center}
where both horizontal sequences are fibre bundles, and we use $O(0^l \oplus\R^{n-l})$ to represent the $O(n-l)$ that comes from the lower right $(n-l)$-by-$(n-l)$ block of $O(n)$ such that the resulted Stiefel manifolds in the lower sequence are the fixed-point sets of the $\Z_2^l$-actions on the upper sequence. 

The following lemma is a mod-$2$ equivariant version of Formula\,(3.9) from Mimura and Toda\,\cite{MT91}.
\begin{lem}\label{lem:mod2EquivCohInduction}
	The pullback $\pi^*:H^*_{\Z^l_2}\big(O(n)/O(k+1),\F_2\big)\rightarrow H^*_{\Z^l_2}\big(O(n)/O(k),\F_2\big)$ is injective. Moreover, $H^*_{\Z^l_2}\big(O(n)/O(k),\F_2\big)$ is a free $H^*_{\Z^l_2}\big(O(n)/O(k+1),\F_2\big)$-module with a basis $\{1,\tilde{h}_{k}\}$ where $\tilde{h}_{k}$ is of degree $k$ and satisfies the fiber integration $\pi_*(\tilde{h}_{k})=1$.
\end{lem}

\begin{proof}
	By Lemma\,\ref{lem:BundleTower}\,(a), the $S^k$-bundle $\pi: O(n)/O(k)\rightarrow O(n)/O(k+1)$ is the unit sphere bundle of the full universal $\R^{k+1}$-bundle $\mathcal{V}\big(V_{n-k}(\R^{n})\big)\rightarrow V_{n-k-1}(\R^{n})$, hence they have the same mod-$2$ $\Z^l_2$-equivariant Euler classes. Recall the complete decomposition\,\eqref{eq:completeDecom}: $V_{n-k-1}(\R^{n})\times \R^n = \big(\oplus_{i=1}^{n-k-1} \mathcal{L}_i\big) \oplus \mathcal{V}$, where $\R^n$ is acted by $\Z^l_2$ flipping the first $l$ coordinates and $\mathcal{L}_i$'s are equivariant trivial line bundles. The Whitney product formula then gives $(1+\beta_1)\cdots (1+\beta_l)=1+\tilde{w}_1(\mathcal{V})+\cdots+\tilde{w}_{k+1}(\mathcal{V})$. Since $k+1>l$, then $\tilde{e}(\mathcal{V})=\tilde{w}_{k+1}(\mathcal{V})=0$. The equivariant Gysin exact sequence of the $S^k$-bundle is
	\begin{center}
		\begin{tikzpicture}[descr/.style={fill=white,inner sep=1.5pt}]
		\matrix (m) [
		matrix of math nodes,
		row sep=2em,
		column sep=2.5em,
		text height=1.5ex, text depth=0.25ex
		]
		{ \cdots & H^{i-k-1}_{\Z^l_2}\Big(\frac{O(n)}{O(k+1)},\F_2\Big) & H^i_{\Z^l_2}\Big(\frac{O(n)}{O(k+1)},\F_2\Big) & H^i_{\Z^l_2}\Big(\frac{O(n)}{O(k)},\F_2\Big)
			& \\
			& & H^{i-k}_{\Z^l_2}\Big(\frac{O(n)}{O(k+1)},\F_2\Big)
			& H^{i+1}_{\Z^l_2}\Big(\frac{O(n)}{O(k+1)},\F_2\Big) & \cdots \\
		};
		
		\path[overlay,->, font=\scriptsize,>=latex]
		(m-1-1) edge node[auto] {$\pi_*$} (m-1-2)
		(m-1-2) edge node[auto] {$\cdot \tilde{e}(\mathcal{V})$} (m-1-3) 
		(m-1-3) edge node[auto] {$\pi^*$} (m-1-4) 
		(m-1-4) edge[out=355,in=175] node[descr,yshift=0.3ex] {$\pi_*$} (m-2-3)
		(m-2-3) edge node[auto] {$\cdot \tilde{e}(\mathcal{V})$} (m-2-4) 
		(m-2-4) edge node[auto] {$\pi^*$} (m-2-5);
		\end{tikzpicture}
	\end{center}
	which reduces to the short exact sequence
	\[
	0{\longrightarrow} H^i_{\Z^l_2}\Big(\frac{O(n)}{O(k+1)},\F_2\Big)\overset{\pi^*}{\longrightarrow} H^i_{\Z^l_2}\Big(\frac{O(n)}{O(k)},\F_2\Big)\overset{\pi_*}{\longrightarrow}H^{i-k}_{\Z^l_2}\Big(\frac{O(n)}{O(k+1)},\F_2\Big){\longrightarrow} 0
	\]
	and implies the current lemma.
\end{proof}

\begin{proof}[Proof of Theorem\,\ref{thm:mod2EquivRealStiefel}]
	Our proof is based on a double induction on $l$ and $n-k$. When $l=0$, then $O(n-l)/O(k-l)=O(n)/O(k)$; when $n-k=0$, then $O(n)/O(k)=pt$. In both cases, Theorem\,\ref{thm:mod2EquivRealStiefel} is true by triviality. For the general case when $1\leq l\leq k \leq n$, suppose Theorem\,\ref{thm:mod2EquivRealStiefel} is true for the $\Z_2^{l}$-action on $O(n)/O(k+1)$ and the $\Z_2^{l-1}$-action on $O(n-1)/O(k-1)$.
	
	Take the $\Z_2^l$-equivariant cohomology rings on the previous diagram of Stiefel manifolds, we have the arrow-reversed diagram
	\begin{center}
		\begin{tikzpicture}[description/.style={fill=white,inner sep=2pt}]
		\matrix (m) [matrix of math nodes, row sep=2em,
		column sep=3em, text height=1.5ex, text depth=0.25ex]
		{H^*_{\Z^l_2}(S^{k},\F_2)  & H^*_{\Z^l_2}\big(O(n)/O(k),\F_2\big) & H^*_{\Z^l_2}\big(O(n)/O(k+1),\F_2\big) \\
			H^*_{\Z^l_2}(pt,\F_2)\otimes_{\F_2}\frac{\F_2[h_{k-l}]}{\big\langle h_{k-l}^2=0 \big\rangle} & H^*_{\Z^l_2}(pt,\F_2)\otimes_{\F_2}\frac{\F_2[h_{i}\mid_{k-l\leq i \leq n-l-1}]}{\big\langle h_{i}^2=h_{2i} \big\rangle} & H^*_{\Z^l_2}(pt,\F_2)\otimes_{\F_2}\frac{\F_2[h_{i}\mid_{k-l+1\leq i \leq n-l-1}]}{\big\langle h_{i}^2=h_{2i} \big\rangle}.\\};
		\path[->]
		(m-1-2) edge node[above]{$j^*$} (m-1-1)
		(m-1-3) edge node[above]{$\pi^*$} (m-1-2)
		(m-2-2) edge node[above]{${j'}^*$}(m-2-1)
		(m-2-3) edge node[above]{${\pi'}^*$} (m-2-2);
		\path[right hook->]
		(m-1-1) edge node[left]{$\iota^*_{k+1,k,l}$} (m-2-1)
		(m-1-2) edge node[left]{$\iota^*_{n,k,l}$} (m-2-2)
		(m-1-3) edge node[left]{$\iota^*_{n,k+1,l}$} (m-2-3);
		\end{tikzpicture}
	\end{center}
	where the lower sequence of cohomology rings are due to the fact that the corresponding Stiefel submanifolds are the $\Z_2^l$-fixed-point sets. On the bottom left corner, the relation $h_{k-l}^2=0$ is true when $k-l\geq 1$ but has to be replaced by the relation $h_{0}^2=h_{0}$ when $k-l=0$. In either case, the following argument is the same.
	
	We have seen that all the vertical homomorphisms in the above diagram are injective. In the upper sequence, $\pi^*$ is injective by Lemma\,\ref{lem:mod2EquivCohInduction}. In the lower sequence, ${\pi'}^*$ is also injective and the homomorphisms ${\pi'}^*,{j'}^*$ map the $H^*_{\Z^l_2}(pt,\F_2)=\F_2[\beta_1,\ldots,\beta_l]$-algebra generators $h_i$'s in a canonical way.
	
	Theorem\,\ref{thm:mod2EquivRealStiefel} is true for the $\Z_2^{l}$-action on $O(n)/O(k+1)$ such that the rightmost vertical monomorphism $\iota^*_{n,k+1,l}$ maps the $\F_2[\beta_1,\ldots,\beta_l]$-generators $\tilde{h}_i$ of $H^*_{\Z^l_2}\Big(O(n)/O(k+1),\F_2\Big)$ as $$\iota^*_{n,k+1,l}\Big(\sum_{i=k+1}^{n-1} \tilde{h}_i\Big)=\Psi_{k+1}\Big(\prod_{j=1}^{l}(1+\beta_j)\cdot\sum_{i=k-l+1}^{n-l-1}h_i\Big)$$  which will be carried over by the monomorphisms $\pi^*,{\pi'}^*$ to be an image of the middle vertical monomorphism: 
	\begin{align}\label{eq:k+1,n-1}
	\iota^*_{n,k,l}\Big(\sum_{i=k+1}^{n-1} \tilde{h}_i\Big)=\Psi_{k+1}\Big(\prod_{j=1}^{l}(1+\beta_j)\cdot\sum_{i=k-l+1}^{n-l-1}h_i\Big)=\Psi_{k+1}\Big(\prod_{j=1}^{l}(1+\beta_j)\cdot\sum_{i=k-l}^{n-l-1}h_i\Big)
	\end{align}
	where the second identity is true because $\prod_{j=1}^{l}(1+\beta_j)\cdot h_{k-l}$ is of degree $k$.
	
	By Lemma\,\ref{lem:mod2EquivCohInduction}, we have the fiber integration $\pi_*(\tilde{h}_{k})=1$, which means that $\tilde{h}_{k}$, if restricted on a sphere fibre $S^{k}$, will be an equivariant unit volume form, i.e. $j^*(\tilde{h}_{k})$ is an $\F_2[\beta_1,\ldots,\beta_l]$-algebra generator of $H^*_{\Z^l_2}(S^{k},\F_2)$. Hence, by Lemma\,\ref{prop:mod2LocalizeSphere}, we can make $\iota^*_{k+1,k,1}(j^*(\tilde{h}_{k}))=(\beta_1\cdots \beta_l)\cdot h_{k-l}$. Going around the left-hand-side square diagram, we have 
	\begin{align}\label{eq:j'k}
	{j'}^*(\iota^*_{n,k,l}(\tilde{h}_{k}))=\iota^*_{k+1,k,l}(j^*(\tilde{h}_{k}))=(\beta_1\cdots \beta_l)\cdot h_{k-l}.
	\end{align}
	
	Consider the $2$-subtori $\Z_2\times \{1\}^{l-1},\,\{1\}\times \Z_2^{l-1} \subset \Z_2^l$, which will be denoted simply as $\Z_2,\Z_2^{l-1}$. The $\Z_2$ acts on $O(n)/O(k)$ with the fixed-point set $O(n-1)/O(k-1)$ on which the $\Z_2^{l-1}$ acts with the fixed-point set $O(n-l)/O(k-l)$. Hence we have the sequence of inclusions of fixed-point sets
	\begin{center}
		\begin{tikzpicture}[description/.style={fill=white,inner sep=2pt}]
		\matrix (m) [matrix of math nodes, row sep=2em,
		column sep=6em, text height=1.5ex, text depth=0.25ex]
		{O(n-l)/O(k-l)  & O(n-1)/O(k-1) & O(n)/O(k) \\};
		\path[right hook->]
		(m-1-1) edge node[above]{$\iota_{n-1,k-1,l-1}$} (m-1-2)
		(m-1-2) edge node[above]{$\iota_{n,k,1}$} (m-1-3);
		\end{tikzpicture}
	\end{center}
	whose composition is the inclusion $\iota_{n,k,l}: O(n-l)/O(k-l) \hookrightarrow O(n)/O(k) $.

	Theorem\,\ref{thm:mod2EquivRealStiefel} is true for the $\Z_2^{l-1}$-action on $O(n-1)/O(k-1)$ such that the $H^*_{\Z_2^{l-1}}(pt,\F_2)=\F_2[\beta_2,\ldots,\beta_l]$-monomorphism $\iota^*_{n-1,k-1,l-1}:H^*_{\Z_{2}^{l-1}}\big(O(n-1)/O(k-1),\F_2\big)\hookrightarrow \F_2[\beta_2,\ldots,\beta_l]\otimes_{\F_2} H^*\big(O(n-l)/O(k-l),\F_2\big)$ maps the $\F_2[\beta_2,\ldots,\beta_l]$-generators $\bar{h}_i$ of $H^*_{\Z^{l-1}_2}\big(O(n-1)/O(k-1),\F_2\big)$ as $$\iota^*_{n-1,k-1,l-1}\Big(\sum_{i=k-1}^{n-2} \bar{h}_i\Big)=\Psi_{k-1}\Big(\prod_{j=2}^{l}(1+\beta_j)\cdot\sum_{i=k-l}^{n-l-1}h_i\Big)$$
	from which we have the components of degrees $k-1,k$
	\begin{align}\label{eq:bar,k-1,k}
	\begin{split}
		\iota^*_{n-1,k-1,l-1}(\bar{h}_{k-1})&=(\beta_2\cdots\beta_l)\cdot h_{k-l}+\sum_{i=1}e_{l-i-1}(\beta_2,\ldots,\beta_l)\cdot h_{k-l+i}\\
	\iota^*_{n-1,k-1,l-1}(\bar{h}_{k})&=\sum_{i=1}e_{l-i}(\beta_2,\ldots,\beta_l)\cdot h_{k-l+i}
	\end{split}
	\end{align}
	where $e_{l-i}(\beta_2,\ldots,\beta_l)$ is the $(l-i)$-th elementary symmetric function in $\beta_2,\ldots,\beta_l$ for $1\leq i \leq l$ and $h_{k-l+i}$ is defined for $k-l\leq k-l+i \leq n-l-1$, the other undefined terms are zero. Since the rank-$1$ $2$-subtorus $\Z_2$ acts trivially on $O(n-1)/O(k-1)$, we have the ring isomorphisms $H^*_{\Z_{2}^l}\big(O(n-1)/O(k-1),\F_2\big)=H^*_{\Z_{2}}(pt,\F_2)\otimes_{\F_2} H^*_{\Z_{2}^{l-1}}\big(O(n-1)/O(k-1),\F_2\big)=\F_2[\beta_1]\otimes_{\F_2}H^*_{\Z_{2}^{l-1}}\big(O(n-1)/O(k-1),\F_2\big)$. Hence, the $\F_2[\beta_2,\ldots,\beta_l]$-generators $\iota^*_{n-1,k-1,l-1}(\bar{h}_i)$ where $k-1\leq i \leq n-2$ are also the $\F_2[\beta_1,\ldots,\beta_l]$-generators for the image of the $\F_2[\beta_1,\ldots,\beta_l]$-monomorphism $ H^*_{\Z_{2}^l}\big(O(n-1)/O(k-1),\F_2\big)\hookrightarrow H^*_{\Z_{2}^l}\big(O(n-l)/O(k-l),\F_2\big)$, which will still be denoted by $\iota^*_{n-1,k-1,l-1}$.
	
	Consider the pullback $\iota_{n,k,1}^*:H^*_{\Z_{2}^l}\big(O(n)/O(k),\F_2\big)\hookrightarrow H^*_{\Z_{2}^l}\big(O(n-1)/O(k-1),\F_2\big)$. By degree reason, the image $\iota_{n,k,1}^*(h_k)$ is an $\F_2$-linear combination of $\bar{h}_k,f\cdot\bar{h}_{k-1},g$ where $f,g\in \F_2[\beta_1,\ldots,\beta_l]$ are homogeneous polynomials of degrees $1,k$. At non-equivariant level, we have $\iota_{n,k,1}^*(h_k)=h_k$, then the equivariant image $\iota^*_{n,k,1}(\tilde{h}_{k})$ must contain an $\bar{h}_k$ term: $\iota^*_{n,k,1}(\tilde{h}_{k})=\bar{h}_k+f\cdot\bar{h}_{k-1}+g$. Applying $\iota^*_{n-1,k-1,l-1}$ on both sides and using the composition $\iota_{n,k,l}=\iota_{n,k,1}\circ\iota_{n-1,k-1,l-1}$, we have
	\begin{align}\label{eq:tilde,bar}
	\iota^*_{n,k,l}(\tilde{h}_{k})=\iota^*_{n-1,k-1,l-1}(\bar{h}_k)+f\cdot\iota^*_{n-1,k-1,l-1}(\bar{h}_{k-1})+g
	\end{align}
	Note ${j'}^*$ annihilates any term involving ${h_{k-l+i},i\geq 1}$, we use \eqref{eq:bar,k-1,k} and get
	\[
	{j'}^*(\iota^*_{n,k,l}(\tilde{h}_{k}))={j'}^*(\iota^*_{n-1,k-1,l-1}(\bar{h}_k))+{j'}^*(f\cdot\iota^*_{n-1,k-1,l-1}(\bar{h}_{k-1}))+{j'}^*(g)=f\cdot (\beta_2\cdots \beta_l) \cdot h_{k-l} + g
	\]
	which, as we have shown in \eqref{eq:j'k}, is supposed to be $(\beta_1\cdots \beta_l) \cdot h_{k-l}$. Hence $f=\beta_1,g=0$. Plug them into \eqref{eq:tilde,bar} and use \eqref{eq:bar,k-1,k}, we have
	\[
	\iota^*_{n,k,l}(\tilde{h}_{k})=\iota^*_{n-1,k-1,l-1}(\bar{h}_k)+\beta_1\cdot\iota^*_{n-1,k-1,l-1}(\bar{h}_{k-1})=(\beta_1\cdots\beta_l)\cdot h_{k-l}+\sum_{i=1}e_{l-i}(\beta_1,\ldots,\beta_l)\cdot h_{k-l+i}
	\]
	which is the degree-$k$ component of $\prod_{j=1}^{l}(1+\beta_j)\cdot\sum_{i=k-l}^{n-l-1}h_i$. Combining with \eqref{eq:k+1,n-1}, we then have verified the formula of $\iota^*_{n,k,l}\big(\sum_{i=k}^{n-1} \tilde{h}_i\big)$.
\end{proof}
	
\vskip 20pt	
\section{$\Z_{1/2}$-coefficient equivariant cohomology rings of the oriented or real flag manifolds}
\vskip 15pt
Given a sequence of positive integers $n_1,\ldots,n_k$, if any two of them, say $n_i,n_j$, are odd, then $\mathrm{rank}\,SO(n_i)+\mathrm{rank}\,SO(n_j)=\mathrm{rank}\,SO(n_i+n_j)-1$. Assuming $n\geq n_1+\cdots+n_k$, the rank difference between $SO(n)$ and $SO(n_1)\times\cdots\times SO(n_k)$ depends on the difference $n-(n_1+\cdots+n_k)$ and also on the parities of all these integers. Because of this complication, in Chap.\,XI of the book of Greub, Halperin and Vanstone\,\cite{GHV76}, the rational ordinary cohomology of some but not all cases of the oriented flag manifolds are considered. We will apply the Leray-Borel Theorem\,\ref{thm:BundleLerayBorel} in several steps to describe the $\Z_{1/2}$-coefficient equivariant cohomology rings of the torus actions on the oriented flag manifolds, then use a covering-space argument to describe the cases of real flag manifolds.

\subsection{Covering between the oriented and real flag manifolds}
The real flag manifold can be identified as 
$
O(n)/\big(O(n_1)\times \cdots \times O(n_k)\big) = SO(n)/S\big(O(n_1)\times \cdots \times O(n_k)\big).
$
Hence we have the covering map between an oriented flag manifold and a real flag manifold:
\[
\frac{SO(n)}{SO(n_1)\times \cdots \times SO(n_k)} \longrightarrow \frac{SO(n)}{S\Big(O(n_1)\times \cdots \times O(n_k)\Big)}
\]
with the covering group $S\big(O(1)^k\big)=\{(\epsilon_1,\ldots,\epsilon_k)\mid \epsilon_j=\pm 1,\, \prod_{j=1}^{k}\epsilon_j=1\}$, which is isomorphic to $\Z_2^{k-1}$. 

Using the notations of Subsubsection\,\ref{subsubsec:generalFlag}, let $(v_1,\ldots,v_{n-N};W_1,\ldots,W_{k}) \in SO(n)/SO(n_1)\times \cdots \times SO(n_k)$. Then the $S\big(O(1)^k\big)$ action is
\[
(\epsilon_1,\ldots,\epsilon_k) \Cdot (v_1,\ldots,v_{n-N};W_1,\ldots,W_{k}) = (v_1,\ldots,v_{n-N};\epsilon_1 W_1,\ldots,\epsilon_1W_{k})
\]
where $\epsilon_j W_j=\pm W_j$ is the space $W_j$ with the original or the reserved orientation. 

The $S\big(O(1)^k\big)$ action on $SO(n)/\big(SO(n_1)\times \cdots \times SO(n_k)\big)$ naturally extends to its universal oriented vector bundles $\mathcal{W}_j$ and commutes with the canonical torus actions, hence induces an action on the equivariant cohomology ring $H_{T^M}^*\big(SO(n)/\big(SO(n_1)\times \cdots \times SO(n_k)\big),{\Z_{1/2}}\big)$. Let $\tilde{p}^{(j)}=1+\tilde{p}_1^{(j)}+\cdots+\tilde{p}_{m_j}^{(j)}$, $\tilde{e}^{(j)}$ be the equivariant total Pontryagin class and the equivariant Euler class of $\mathcal{W}_j$. Since Pontryagin classes ignore orientations of the vector bundles, while Euler classes detect the orientations, we have
\begin{align*}
(\epsilon_1,\ldots,\epsilon_k) \Cdot \tilde{p}^{(j)}=\tilde{p}^{(j)} \qquad \qquad (\epsilon_1,\ldots,\epsilon_k) \Cdot \tilde{e}^{(j)}=\epsilon_j \tilde{e}^{(j)}.
\end{align*}

We need the following well-known fact about the cohomology of covering spaces.
\begin{lem}\label{lem:coveringSpace}
	Let $\pi:X\rightarrow Y$ be a covering between compact topological spaces with a finite covering group $\Gamma$. There is an induced $\Gamma$-action on the cohomology $H^*(X,\K)$. If $|\Gamma|$ is invertible in $\K$, then $\pi^*: H^*(Y,\K) \rightarrow H^*(X,\K)$ is injective and its image is the $\Gamma$-invariant subring $H^*(X,\K)^\Gamma\subset H^*(X,\K)$. This statement is also true for equivariant cohomology if a torus $T$ acts on $X$ and commutes with the action of $\Gamma$. 
\end{lem}

Since $|S\big(O(1)^k\big)|=2^{k-1}$ is invertible in $\Z_{1/2}$, in the following discussions, as soon as we derive the $\Z_{1/2}$-coefficient equivariant cohomology rings of the oriented flag manifolds, we then derive the equivariant cohomology rings of the real flag manifolds by finding the $S\big(O(1)^k\big)$-invariants.

\subsection{The case of real Stiefel manifolds}
By Proposition\,\ref{prop:FormalStiefel}, the $T^l$-actions on the real Stiefel manifolds $SO(2m)/SO(2k+1),SO(2m+1)/SO(2k+1),SO(2m)/SO(2k),SO(2m+1)/SO(2k)$, where $l\leq k$, are equivariantly formal in $\Z_{1/2}$ coefficients. Together using Theorem\,\ref{thm:realStief}\,(b), we have $H_{T^l}^*(pt,\Z_{1/2})=\Z_{1/2}[\alpha_1,\ldots,\alpha_l]$-module isomorphisms
\begin{align*}
H_{T^{l}}^*(SO(2m)/SO(2k+1),{\Z_{1/2}})&\cong {\Z_{1/2}}[\alpha_1,\ldots,\alpha_{l}]\otimes_{\Z_{1/2}}\Lambda_{\Z_{1/2}}[y_{k+1},\ldots,y_{m-1},x_m]\\
H_{T^{l}}^*(SO(2m+1)/SO(2k+1),{\Z_{1/2}})&\cong {\Z_{1/2}}[\alpha_1,\ldots,\alpha_{l}]\otimes_{\Z_{1/2}}\Lambda_{\Z_{1/2}}[y_{k+1},\ldots,y_{m}]\\
H_{T^{l}}^*(SO(2m)/SO(2k),{\Z_{1/2}})&\cong {\Z_{1/2}}[\alpha_1,\ldots,\alpha_{l}]\otimes_{\Z_{1/2}}\Lambda_{\Z_{1/2}}[e_k, y_{k+1},\ldots,y_{m-1},x_m]\\
H_{T^{l}}^*(SO(2m+1)/SO(2k),{\Z_{1/2}})&\cong {\Z_{1/2}}[\alpha_1,\ldots,\alpha_{l}]\otimes_{\Z_{1/2}}\Lambda_{\Z_{1/2}}[e_k,y_{k+1},\ldots,y_{m}].
\end{align*}

We can improve the $\Z_{1/2}[\alpha_1,\ldots,\alpha_l]$-module isomorphisms to be $\Z_{1/2}[\alpha_1,\ldots,\alpha_l]$-algebra isomorphisms. In the following theorem, we consider the case when $l=k$.
\begin{thm}\label{thm:equivRealStiefel}
	For the canonical $T^k$-actions on the real Stiefel manifolds $SO(2m)/SO(2k+1),SO(2m+1)/SO(2k+1),SO(2m)/SO(2k),SO(2m+1)/SO(2k)$, 
	\begin{enumerate}[label=(\alph*)]
		\item as $\Z_{1/2}[\alpha_1,\ldots,\alpha_{k}]$-algebras,
		\begin{align*}
		H_{T^{k}}^*(SO(2m)/SO(2k+1),{\Z_{1/2}})&= {\Z_{1/2}}[\alpha_1,\ldots,\alpha_{k}]\otimes_{\Z_{1/2}}\Lambda_{\Z_{1/2}}[\tilde{y}_{k+1},\ldots,\tilde{y}_{m-1},\tilde{x}_m]\\
		H_{T^{k}}^*(SO(2m+1)/SO(2k+1),{\Z_{1/2}})&= {\Z_{1/2}}[\alpha_1,\ldots,\alpha_{k}]\otimes_{\Z_{1/2}}\Lambda_{\Z_{1/2}}[\tilde{y}_{k+1},\ldots,\tilde{y}_{m}]
		\end{align*}
		where $\tilde{x}_m,\tilde{y}_i$ are of odd degrees $2m-1,4i-1$;
		\item the full universal vector bundles $\mathcal{V}(SO(2m)/SO(2k))$ and $\mathcal{V}(SO(2m+1)/SO(2k))$ are of rank $2k$ and denote both of the $T^{k}$-equivariant Euler classes by $\tilde{e}_k$, then it satisfies the quadratic equation
		\[
		\tilde{e}_k^2=\alpha_1^2\cdots \alpha_{k}^2;
		\]
		\item
		as ${\Z_{1/2}}[\alpha_1,\ldots,\alpha_{k}]$-algebras 
		\begin{align*}
		H_{T^{k}}^*(SO(2m)/SO(2k),{\Z_{1/2}})
		&= \frac{{\Z_{1/2}}[\alpha_1,\ldots,\alpha_{k};\tilde{e}_k]}{\big\langle \tilde{e}_k^2=\alpha_1^2\cdots \alpha_{k}^2 \big\rangle} \otimes_{\Z_{1/2}}\Lambda_{\Z_{1/2}}[\tilde{y}_{k+1},\ldots,\tilde{y}_{m-1},\tilde{x}_m]\\
		H_{T^{k}}^*(SO(2m+1)/SO(2k),{\Z_{1/2}})
		&= \frac{{\Z_{1/2}}[\alpha_1,\ldots,\alpha_{k};\tilde{e}_k]}{\big\langle \tilde{e}_k^2=\alpha_1^2\cdots \alpha_{k}^2 \big\rangle}\otimes_{\Z_{1/2}}\Lambda_{\Z_{1/2}}[\tilde{y}_{k+1},\ldots,\tilde{y}_{m}].
		\end{align*}
	\end{enumerate}  
\end{thm}
\begin{proof}
	Via the ${\Z_{1/2}}[\alpha_1,\ldots,\alpha_{k}]$-module isomorphisms, we can lift the generators $x_{m},y_i,e_k$ of the ordinary cohomology rings to be ${\Z_{1/2}}[\alpha_1,\ldots,\alpha_{k}]$-generators $\tilde{x}_{m},\tilde{y}_i,\tilde{e}_k$ of the equivariant cohomology rings.
	\begin{enumerate}[label=(\alph*)]
		\item For the odd-degree generators $\tilde{x}_{m},\tilde{y}_i$, we get $\tilde{x}_{m}^2=-\tilde{x}_{m}^2,\,\tilde{y}_i^2=-\tilde{y}_i^2$. Since we are working with $\Z_{1/2}$ coefficients, there are no $2$-torsions, hence we have $\tilde{x}_{m}^2=0,\,\tilde{y}_i^2=0$. 
		\item Using one of Relations\,\eqref{eq:EquivWhitney}, we have the $T^k$-equivariant total Pontryagin classes of the rank-$2k$ universal bundles $\mathcal{V}(SO(2m)/SO(2k)),\mathcal{V}(SO(2m+1)/SO(2k))$:
		\[
		1+\tilde{p}_1+\cdots+\tilde{p}_k = (1+\alpha_1^2)\cdots(1+\alpha_k^2)
		\]
		from which we get $\tilde{p}_k=\alpha_1^2\cdots \alpha_k^2$. The square formula of Euler class then gives
		$
		\tilde{e}_k^2=\tilde{p}_k=\alpha_1^2\cdots \alpha_k^2.
		$
		\item By (a), we already have $\tilde{x}_{m}^2=0,\,\tilde{y}_i^2=0$. By (b), we just need to update the quadratic equation $e_k^2=0$ to $\tilde{e}_k^2=\alpha_1^2\cdots \alpha_k^2$. 	
	\end{enumerate}
\end{proof}

\begin{rmk}\label{rmk:equivRealStiefel1}
	Comparing the above (a) and (c), we get a ring isomorphism
	\begin{align*}
	H_{T^{k}}^*(SO(n)/SO(2k),{\Z_{1/2}}) = \frac{H_{T^{k}}^*(SO(n)/SO(2k+1),{\Z_{1/2}})[\tilde{e}_{k}]}{\big\langle \tilde{e}_k^2=\alpha_1^2\cdots \alpha_k^2\big\rangle}.
	\end{align*}
\end{rmk}

\begin{rmk}\label{rmk:equivRealStiefel2}
	Because of equivariant formality, if $l<k$, we simply substitute $\alpha_{l+1}=\cdots=\alpha_k=0$ into the above Theorem\,\ref{thm:equivRealStiefel} and get the $\Z_{1/2}[\alpha_1,\ldots,\alpha_l]$-algebras
	\begin{align*}
	H_{T^{l}}^*(SO(2m)/SO(2k+1),{\Z_{1/2}})&= {\Z_{1/2}}[\alpha_1,\ldots,\alpha_{l}]\otimes_{\Z_{1/2}}\Lambda_{\Z_{1/2}}[\tilde{y}_{k+1},\ldots,\tilde{y}_{m-1},\tilde{x}_m]\\
	H_{T^{l}}^*(SO(2m+1)/SO(2k+1),{\Z_{1/2}})&= {\Z_{1/2}}[\alpha_1,\ldots,\alpha_{l}]\otimes_{\Z_{1/2}}\Lambda_{\Z_{1/2}}[\tilde{y}_{k+1},\ldots,\tilde{y}_{m}]\\
	H_{T^{l}}^*(SO(2m)/SO(2k),{\Z_{1/2}})&= {\Z_{1/2}}[\alpha_1,\ldots,\alpha_{l}]\otimes_{\Z_{1/2}}\Lambda_{\Z_{1/2}}[\tilde{e}_k, \tilde{y}_{k+1},\ldots,\tilde{y}_{m-1},\tilde{x}_m]\\
	H_{T^{l}}^*(SO(2m+1)/SO(2k),{\Z_{1/2}})&= {\Z_{1/2}}[\alpha_1,\ldots,\alpha_{l}]\otimes_{\Z_{1/2}}\Lambda_{\Z_{1/2}}[\tilde{e}_k,\tilde{y}_{k+1},\ldots,\tilde{y}_{m}].
	\end{align*}
\end{rmk}

Similar to the localization results (Theorem\,\ref{thm:mod2EquivRealStiefel}, Proposition\,\ref{prop:mod2LocalizeSphere}) of the $2$-torus actions on real Stiefel manifolds, we have the following localization results for the torus actions on real Stiefel manifolds. The proofs are technically the same, so we omit the details.

\begin{prop}\label{prop:localizeSphere}
	The $T^l$-actions on $SO(2k+1)/SO(2k)=S^{2k},\,SO(2k+2)/SO(2k+1)=S^{2k+1},\,SO(2k+3)/SO(2k+1)$, where $l\leq k$, have the fixed-point sets $S^{2k-2l},\,S^{2k-2l+1},\,SO(2k-2l+3)/SO(2k-2l+1)$. The inclusions of the $T^l$-fixed-point sets induce the pullback monomorphisms 
	\begin{align*}
	H^*_{T^l}(S^{2k},\Z_{1/2})&\longhookrightarrow H^*_{T^l}(S^{2k-2l},\Z_{1/2})=
	\begin{cases}
	\Z_{1/2}[\alpha_1,\ldots,\alpha_l][e_{k-l}]/\big\langle e_{k-l}^2\big\rangle, & k-l\geq 1\\
	\Z_{1/2}[\alpha_1,\ldots,\alpha_l][e_0]/\big\langle e_0^2-e_0\big\rangle, & k-l=0,
	\end{cases}\\
	H^*_{T^l}(S^{2k+1},\Z_{1/2})&\longhookrightarrow H^*_{T^l}(S^{2k-2l+1},\Z_{1/2})=
	\Z_{1/2}[\alpha_1,\ldots,\alpha_l]\otimes_{\Z_{1/2}}\Lambda_{\Z_{1/2}}[x_{k-l+1}]\\
	H^*_{T^l}\Big(\frac{SO(2k+3)}{SO(2k+1)},\Z_{1/2}\Big)&\longhookrightarrow H^*_{T^l}\Big(\frac{SO(2k-2l+3)}{SO(2k-2l+1)},\Z_{1/2}\Big)=
	\Z_{1/2}[\alpha_1,\ldots,\alpha_l]\otimes_{\Z_{1/2}}\Lambda_{\Z_{1/2}}[y_{k-l+1}]
	\end{align*}
	whose images, as $\Z_{1/2}[\alpha_1,\ldots,\alpha_l]$-subalgebras, are generated on $(\alpha_1\cdots\alpha_l) \cdot e_{k-l},\,(\alpha_1\cdots\alpha_l) \cdot x_{k-l+1},\,(\alpha_1\cdots\alpha_l)^2 \cdot y_{k-l+1}$ respectively. 
\end{prop}

\begin{thm}\label{thm:localizeRealStiefel}
	The $T^l$-action on $SO(2m)/SO(2k)$, where $l\leq k$, has the fixed-point set $SO(2m-2l)/SO(2k-2l)$. The inclusion $\iota:SO(2m-2l)/SO(2k-2l) \hookrightarrow SO(2m)/SO(2k)$ induces the pullback monomorphism
	\[
	\iota^*:H^*_{T^l}\Big(\frac{SO(2m)}{SO(2k)},\Z_{1/2}\Big) \longhookrightarrow H^*_{T^l}\Big(\frac{SO(2m-2l)}{SO(2k-2l)},\Z_{1/2}\Big)=\Z_{1/2}[\alpha_1,\ldots,\alpha_l]\otimes_{\Z_{1/2}}\Lambda_{\Z_{1/2}}[y_{k-l+1},\ldots,y_{m-l-1},x_{m-l}]
	\]
	such that the $\Z_{1/2}[\alpha_1,\ldots,\alpha_l]$-generators of $H^*_{T^l}(SO(2m)/SO(2k),\Z_{1/2})$ can be identified by the formulas:
	\begin{align*}
	\iota^*(\tilde{e}_k)&=(\alpha_1\cdots\alpha_l) \cdot e_{k-l}\\
	\iota^*(\tilde{x}_m)&=(\alpha_1\cdots\alpha_l) \cdot x_{m-l}\\
	\iota^*\Big(\sum_{i=k+1}^{n-1} \tilde{y}_i\Big)&=\Psi_{4k+3}\Big(\prod_{j=1}^{l}(1+\alpha^2_j)\cdot\sum_{i=k-l+1}^{n-l-1}y_i\Big).
	\end{align*}
	The cases of $T^l$-actions on $SO(2m)/SO(2k+1),SO(2m+1)/SO(2k+1),SO(2m+1)/SO(2k)$ are similar.
\end{thm}

\subsection{The case when $n_1,\ldots,n_k$ are all even}
In this subsection, suppose $(n_1,\ldots,n_k)=(2m_1,\ldots,2m_k)$ are all even and set $M=m_1+\cdots +m_k$. 

\begin{thm}\label{thm:orientFlag1}
	Consider the canonical $T^M$-action on the oriented flag manifold $SO(n)/\big(SO(2m_1)\times\cdots\times SO(2m_k)\big)$. As an $H^*_{T^M}(pt,{\Z_{1/2}})={\Z_{1/2}}[\alpha_1,\ldots,\alpha_M]$-algebra,
	\begin{enumerate}[label=(\alph*)]
		\item if $n=2M$, then $H_{T^M}^*\big({SO(n)}/\big({SO(2m_1)\times \cdots\times SO(2m_k)}\big),{\Z_{1/2}}\big)$ is
		\begin{align*}
		\frac{{\Z_{1/2}}[\alpha_l\mid_{1\leq l\leq M}][\tilde{p}_i^{(j)}\mid_{1\leq j\leq k,1\leq i\leq m_j}][\tilde{e}^{(j)} \mid_{1\leq j\leq k}]}{\big\langle\prod_{j=1}^{k}\tilde{p}^{(j)}=\prod_{l=1}^{M}(1+\alpha^2_l); (\tilde{e}^{(j)})^2=\tilde{p}_{m_j}^{(j)} \mid_{1\leq j\leq k};\prod_{j=1}^k \tilde{e}^{(j)}=\prod_{l=1}^{M}\alpha_l \big\rangle}
		\end{align*}
		where $\tilde{p}^{(j)}=1+\tilde{p}_1^{(j)}+\cdots+\tilde{p}_{m_j}^{(j)}$, $\tilde{e}^{(j)}$ are the equivariant total Pontryagin class and the equivariant Euler class of the $j$-th universal bundle over $SO(n)/\big(SO(2m_1)\times\cdots\times SO(2m_k)\big)$;
		\item if $n>2M$, then $H_{T^M}^*\big({SO(n)}/\big({SO(2m_1)\times \cdots\times SO(2m_k)}\big),{\Z_{1/2}}\big)$ is
		\begin{align*}
		\frac{{\Z_{1/2}}[\alpha_l\mid_{1\leq l\leq M}][\tilde{p}_i^{(j)}\mid_{1\leq j\leq k,1\leq i\leq m_j}][\tilde{e}^{(j)} \mid_{1\leq j\leq k}]}{\big\langle\prod_{j=1}^{k}\tilde{p}^{(j)}=\prod_{l=1}^{M}(1+\alpha^2_l); (\tilde{e}^{(j)})^2=\tilde{p}_{m_j}^{(j)} \mid_{1\leq j\leq k} \big\rangle} \otimes_{\Z_{1/2}} H^*\big( \frac{SO(n)}{SO(2M+1)},{\Z_{1/2}}\big)
		\end{align*}
		where $H^*( {SO(n)}/{SO(2M+1)},{\Z_{1/2}})$ is an exterior algebra generated by odd-degree elements as in Theorem\,\ref{thm:realStief}\,(a).
	\end{enumerate}
\end{thm}

\begin{proof}
	Since $SO(2M)$ and its subgroup $SO(2m_1)\times \cdots\times SO(2m_k)$ are of the same rank $M$, we use Remark\,\ref{rmk:BundleLerayBorel}\,(b) to get the cohomology ring $H^*\big({SO(2M)}/\big({SO(2m_1)\times \cdots\times SO(2m_k)}\big),{\Z_{1/2}}\big)$ as
	\begin{align*}
	\frac{H^*\big(BSO(2m_1),{\Z_{1/2}}\big)\otimes_{\Z_{1/2}}\cdots \otimes_{\Z_{1/2}} H^*\big(BSO(2m_k),{\Z_{1/2}}\big)}{\big\langle H^+(BSO(2M),{\Z_{1/2}}) \big\rangle}&=\frac{{\Z_{1/2}}[{p}_i^{(j)}\mid_{1\leq j\leq k,1\leq i\leq m_j}][{e}^{(j)} \mid_{1\leq j\leq k}]}{\big\langle\prod_{j=1}^{k}{p}^{(j)}=1; ({e}^{(j)})^2={p}_{m_j}^{(j)} \mid_{1\leq j\leq k};\prod_{j=1}^k {e}^{(j)}=0 \big\rangle}
	\end{align*}
	where the relations in the denominator are exactly Relations\,\eqref{eq:Whitney},\,\eqref{eq:sqEuler},\eqref{eq:WhitneyEuler} for ${SO(2M)}/\big({SO(2m_1)\times \cdots\times SO(2m_k)}\big)$.
	\begin{enumerate}[label=(\alph*)]
		\item If $n=2M$, since $H^*\big({SO(2M)}/\big({SO(2m_1)\times \cdots\times SO(2m_k)}\big),{\Z_{1/2}}\big)$ has only even-degree elements in the above expression, the $T^M$-action on ${SO(2M)}/\big({SO(2m_1)\times \cdots\times SO(2m_k)}\big)$ is equivariantly formal by Remark\,\ref{rmk:formality}. As a $\Z_{1/2}[\alpha_1,\ldots,\alpha_M]$-module, $H_{T^M}^*\big({SO(2M)}/\big({SO(2m_1)\times \cdots\times SO(2m_k)}\big),{\Z_{1/2}}\big)$ is isomorphic to
		\[
		{\Z_{1/2}}[\alpha_1,\ldots,\alpha_M]\otimes_{\Z_{1/2}}H^*\Big(\frac{SO(2M)}{SO(2m_1)\times \cdots\times SO(2m_k)},{\Z_{1/2}}\Big)
		\]
		To get a ${\Z_{1/2}}[\alpha_1,\ldots,\alpha_M]$-algebra isomorphism, we lift the characteristic classes ${p}_i^{(j)},e^{(j)}$ and their relations to be the equivariant characteristic classes and the equivariant Relations\,\eqref{eq:EquivWhitney},\,\eqref{eq:EquivsqEuler}, and Relation\,\eqref{eq:EquivWhitneyEuler}: $\prod_{j=1}^k\tilde{e}^{(j)}=e^{T^M}(\tilde{\R}^{2M})=\prod_{l=1}^{M}\alpha_l$.
		\item If $n>2M$, consider the $T^M$-equivariant bundle 
		\[
		\frac{SO(2M)}{SO(2m_1)\times \cdots \times SO(2m_k)}
		\overset{\iota}{\longhookrightarrow}
		\frac{SO(n)}{SO(2m_1)\times \cdots \times SO(2m_k)}
		\overset{\pi}{\longrightarrow}
		\frac{SO(n)}{SO(2M)}.
		\]
		By Remark\,\ref{rmk:BundleLerayBorel}\,(c), we can use an equivariant version of the Leray-Borel Theorem\,\ref{thm:BundleLerayBorel} to express the ${\Z_{1/2}}[\alpha_1,\ldots,\alpha_M]$-module structure of $H_{T^M}^*\big({SO(n)}/\big({SO(2m_1)\times \cdots\times SO(2m_k)\big)},{\Z_{1/2}}\big)$ as 
		\begin{align*}
		 H_{T^M}^*\Big(\frac{SO(n)}{SO(2M)},{\Z_{1/2}}\Big)\otimes_{\Z_{1/2}}H^*\Big(\frac{SO(2M)}{SO(2m_1)\times \cdots\times SO(2m_k)},{\Z_{1/2}}\Big)
		\end{align*}
		from which we see the $T^M$-action on ${SO(n)}/\big({SO(2m_1)\times \cdots\times SO(2m_k)}\big)$ is equivariantly formal. To get a ${\Z_{1/2}}[\alpha_1,\ldots,\alpha_M]$-algebra isomorphism, we use the equivariant Relations\,\eqref{eq:EquivWhitney},\,\eqref{eq:EquivsqEuler}, and the Relation\,\eqref{eq:EquivWhitneyEuler}: $\prod_{j=1}^{k}\tilde{e}^{(j)}=\tilde{e}(\mathcal{V})$ where the full universal bundle $\mathcal{V}$ stands on both ${SO(n)}/\big({SO(2m_1)\times \cdots\times SO(2m_k)}\big)$ and ${SO(n)}/{SO(2M)}$. By Theorem\,\ref{thm:equivRealStiefel}\,(c), the even-degree generator of $H_{T^M}^*\big({SO(n)}/{SO(2M)},{\Z_{1/2}}\big)$ is $\tilde{e}(\mathcal{V})$, and then is a redundant generator in $H_{T^M}^*\big({SO(n)}/\big({SO(2m_1)\times \cdots\times SO(2m_k)}\big),{\Z_{1/2}}\big)$ because $\tilde{e}(\mathcal{V})=\prod_{j=1}^{k}\tilde{e}^{(j)}$. If we clear out $\tilde{e}(\mathcal{V})$ from $H_{T^M}^*\big({SO(n)}/{SO(2M)},{\Z_{1/2}}\big)$, then by Remark\,\ref{rmk:equivRealStiefel1}, the resulted ${\Z_{1/2}}[\alpha_1,\ldots,\alpha_M]$-subalgebra is $H_{T^M}^*\big({SO(n)}/{SO(2M+1)},{\Z_{1/2}}\big)$. Finally, by Theorem\,\ref{thm:equivRealStiefel}\,(a), $H_{T^M}^*\big({SO(n)}/{SO(2M+1)},{\Z_{1/2}}\big)={\Z_{1/2}}[\alpha_1,\ldots,\alpha_M]\otimes_{\Z_{1/2}}H^*\big({SO(n)}/{SO(2M+1)},{\Z_{1/2}}\big)$, then we get an $H^*\big({SO(n)}/{SO(2M+1)},{\Z_{1/2}}\big)$-factor for $H_{T^M}^*\big({SO(n)}/\big({SO(2m_1)\times \cdots\times SO(2m_k)}\big),{\Z_{1/2}}\big)$.
	\end{enumerate}
\end{proof}

\begin{thm}\label{thm:realFlag1}
		Consider the canonical $T^M$-action on the real flag manifold $O(n)/\big(O(2m_1)\times\cdots\times O(2m_k)\big)$. As a ${\Z_{1/2}}[\alpha_1,\ldots,\alpha_M]$-algebra,
	\begin{enumerate}[(a)]
		\item if $n=2M$, then $H_{T^M}^*\big({O(n)}/\big({O(2m_1)\times \cdots\times O(2m_k)}\big),{\Z_{1/2}}\big)$ is
		\begin{align*}
		\frac{{\Z_{1/2}}[\alpha_l\mid_{1\leq l\leq M}][\tilde{p}_i^{(j)}\mid_{1\leq j\leq k,1\leq i\leq m_j}]}{\big\langle\prod_{j=1}^{k}\tilde{p}^{(j)}=\prod_{l=1}^{M}(1+\alpha^2_l) \big\rangle};
		\end{align*}
		\item if $n>2M$,  then $H_{T^M}^*\big({O(n)}/\big({O(2m_1)\times \cdots\times O(2m_k)}\big),{\Z_{1/2}}\big)$ is
		\begin{align*}
		\frac{{\Z_{1/2}}[\alpha_l\mid_{1\leq l\leq M}][\tilde{p}_i^{(j)}\mid_{1\leq j\leq k,1\leq i\leq m_j}][\tilde{e}_M]}{\big\langle\prod_{j=1}^{k}\tilde{p}^{(j)}=\prod_{l=1}^{M}(1+\alpha^2_l);\tilde{e}_M^2=\prod_{l=1}^{M}\alpha^2_l \big\rangle} \otimes_{\Z_{1/2}} H^*\Big( \frac{SO(n)}{SO(2M+1)},{\Z_{1/2}}\Big)
		\end{align*}
		where $\tilde{e}_M$ is the equivariant Euler class of the rank-$2M$ full universal bundle $\mathcal{V}$ which is oriented and pulled back via the map ${SO(n)}/S\big({O(2m_1)\times \cdots\times O(2m_k)}\big)\rightarrow SO(n)/SO(2M)$.
	\end{enumerate}
\end{thm}
\begin{proof}
	We need to find all the $S\big(O(1)^k\big)$-invariants of $H_{T^M}^*\big({SO(n)}/\big({SO(2m_1)\times \cdots\times SO(2m_k)}\big),{\Z_{1/2}}\big)$. First, the equivariant Pontryagin classes are $S\big(O(1)^k\big)$-invariants. Since the squares of Euler classes are the corresponding top Pontryagin classes, it remains to consider the products of distinct equivariant Euler classes and the generators of the exterior algebra.
	
	We have $k$ distinct equivariant Euler classes $\tilde{e}^{(1)},\ldots,\tilde{e}^{(k)}$. For any proper subset $\{i_1,\ldots,i_l\}\subset \{1,\ldots,k\}$ where $l<k$, choose $i_0\in \{1,\ldots,k\}\setminus\{i_1,\ldots,i_l\}$, and set $\epsilon_{i_0}=\epsilon_{i_1}=-1$ while all the other $\epsilon_i=1$, then 
	\[
	(\epsilon_1,\ldots,\epsilon_k) \Cdot \prod_{j=1}^l \tilde{e}^{(i_j)}=-\prod_{j=1}^l \tilde{e}^{(i_j)},
	\]
	i.e. the partial product $\prod_{j=1}^l \tilde{e}^{(i_j)}$ is not an $S\big(O(1)^k\big)$-invariant. However, by definition $\prod_{j=1}^{k}\epsilon_j=1$, then 
	\[
	(\epsilon_1,\ldots,\epsilon_k) \Cdot \prod_{j=1}^{k}\tilde{e}^{(j)} = \prod_{j=1}^{k}\epsilon_j \prod_{j=1}^{k}\tilde{e}^{(j)} = \prod_{j=1}^{k}\tilde{e}^{(j)}
	\]
	i.e. the full product $\prod_{j=1}^{k}\tilde{e}^{(j)}$ is an $S\big(O(1)^k\big)$-invariant.
	\begin{enumerate}[(a)]
		\item Corresponding to case\,(a) of Theorem\,\ref{thm:orientFlag1}, the full product $\prod_{j=1}^{k}\tilde{e}^{(j)}$ can be replaced by $\prod_{l=1}^{M}\alpha_l$. 
		Hence, the $S\big(O(1)^k\big)$-invariants are entirely generated by the equivariant Pontryagin classes.
		\item Corresponding to case\,(b) of Theorem\,\ref{thm:orientFlag1}, by Relation\,\eqref{eq:WhitneyEuler}, the full product $\prod_{i=1}^{k}\tilde{e}^{(j)}$ is the equivariant Euler class $\tilde{e}_M$ of the full universal bundle $\mathcal{V}$ over ${SO(n)}/\big({SO(2m_1)\times \cdots\times SO(2m_k)}\big)$. Since $\mathcal{V}$ is pulled back from ${SO(n)}/{SO(2M)}$, so is the relation $\tilde{e}_M^2=\prod_{l=1}^{M}\alpha^2_l$ as in Theorem\,\ref{thm:equivRealStiefel}\,(b). Hence $\tilde{e}_M$ is nonzero and can be treated as one of the $S\big(O(1)^k\big)$-invariant generators. Note that $S\big(O(1)^k\big) \subset S\big(O(2m_1)\times \cdots \times O(2m_k)\big) \subset SO(2M) \subset SO(2M+1)$, the action of $S\big(O(1)^k\big)$ on ${SO(n)}/{SO(2M+1)}$ is trivial. Hence the exterior algebra $H^*( {SO(n)}/{SO(2M+1)},{\Z_{1/2}})$ is $S\big(O(1)^k\big)$-invariant. 
	\end{enumerate}
\end{proof}

\begin{cor}\label{cor:flag1}
	The $\Z_{1/2}$-coefficient ordinary cohomology rings of the oriented and flag manifolds in Theorems\,\ref{thm:orientFlag1},\,\ref{thm:realFlag1} can be obtained by substituting $\alpha_1=\cdots=\alpha_M=0$. The ordinary Poincare polynomials of the two oriented flag manifolds in Theorem\,\ref{thm:orientFlag1} are respectively
	\[
	\binom{M}{m_1,\ldots,m_k}_{t^4}\cdot\frac{\prod_{j=1}^k (1+t^{2m_j})}{1+t^{2M}}, \qquad \qquad 
	\binom{M}{m_1,\ldots,m_k}_{t^4}\cdot\Big(\prod_{j=1}^k (1+t^{2m_j})\Big) \cdot P\Big( \frac{SO(n)}{SO(2M+1)},t\Big).
	\]
	The ordinary Poincare polynomials of the two real flag manifolds in Theorem\,\ref{thm:realFlag1} are respectively
	\[
	\binom{M}{m_1,\ldots,m_k}_{t^4}, \qquad \qquad \qquad \qquad \qquad
	\binom{M}{m_1,\ldots,m_k}_{t^4}\cdot (1+t^{2M}) \cdot P\Big( \frac{SO(n)}{SO(2M+1)},t\Big).
	\]
	All these four cases of flag manifolds are orientable. 
\end{cor}
\begin{proof}
	The first sentence on the ordinary cohomology is due to the equivariant formality of the canonical $T^M$-actions. 
	
	To compute the Poincare polynomials, we view the cohomology as built from three steps: the ring of Pontryagin classes subject to the Whitney product formula, the extension by Euler classes subject to the square formula and Whitney product formula, and the exterior algebra. Then the Poincare polynomial is the multiplication of these three contributions.
	
	The oriented flag manifolds are certainly orientable. Since $M=m_1+\cdots+m_k$, the Poincare polynomials of cases\,(a) of both Theorem\,\ref{thm:orientFlag1} and Theorem\,\ref{thm:realFlag1} are of the same degree which is the top dimension of the corresponding real flag manifold, then the non-vanishing of the top-dimensional cohomology implies the orientability. The orientability of case\,(b) of Theorem\,\ref{thm:realFlag1} is similar. 
\end{proof}

\subsection{The case when $n_1,\ldots,n_k$ are not all even}
In this subsection, suppose $(n_1,\ldots,n_k)=(2m_1+1,\ldots,2m_l+1,2m_{l+1},\ldots,2m_k)$ with $l\geq 1$ odd integers. Set $M=m_1+\cdots +m_k$, then $n\geq \sum_i n_i\geq 2M+1$.

Consider the canonical $T^M$-actions on the flag manifolds 
\begin{align*}
X&=\frac{SO(n)}{SO(2m_1+1)\times \cdots \times SO(2m_{l}+1)\times SO(2m_{l+1})\times \cdots \times SO(2m_k)}\\
Y&=\frac{O(n)}{O(2m_1+1)\times \cdots \times O(2m_{l}+1)\times O(2m_{l+1})\times \cdots \times O(2m_k)}.
\end{align*}

\begin{thm}\label{thm:orientFlag2}
As a ${\Z_{1/2}}[\alpha_1,\ldots,\alpha_M]$-algebra, $H_{T^M}^*(X,{\Z_{1/2}})$ is
		\begin{align*}
		\frac{{\Z_{1/2}}[\alpha_r\mid_{1\leq r\leq M}][\tilde{p}_i^{(j)}\mid_{1\leq j\leq k,1\leq i\leq m_j}][\tilde{e}^{(j)} \mid_{l+1\leq j\leq k}]}{\big\langle\prod_{j=1}^{k}\tilde{p}^{(j)}=\prod_{r=1}^{M}(1+\alpha^2_r); (\tilde{e}^{(j)})^2=\tilde{p}_{m_j}^{(j)} \mid_{l+1\leq j\leq k} \big\rangle} \otimes_{\Z_{1/2}} H^*\Big( \frac{SO(n)}{SO(2M+1)},{\Z_{1/2}}\Big).
		\end{align*}
\end{thm}
\begin{proof}
Consider the $T^M$-equivariant bundle projection
		\[
		\frac{SO(n)}{SO(2m_1) \times\cdots\times SO(2m_k)} \longrightarrow
		X
		\]
		whose fibre is 
		\[\frac{SO(2m_1+1)}{SO(2m_1)} \times \cdots \times \frac{SO(2m_l+1)}{SO(2m_l)}\cong S^{2m_1}\times\cdots\times S^{2m_l}\]
		which satisfies the equal-rank condition. By the Leray-Borel Theorem\,\ref{thm:BundleLerayBorel} and Remark\,\ref{rmk:BundleLerayBorel}\,(3), as an $H_{T^M}^*(X,{\Z_{1/2}})$-module, $H_{T^M}^*\big({SO(n)}/\big({SO(2m_1)\times \cdots\times SO(2m_k)}\big),{\Z_{1/2}}\big)$ is isomorphic to
		\begin{align*}
		H_{T^M}^*(X,{\Z_{1/2}}) \otimes_{\Z_{1/2}} H^*(S^{2m_1}\times\cdots\times S^{2m_l},{\Z_{1/2}}) \cong H_{T^M}^*(X,{\Z_{1/2}}) \otimes_{\Z_{1/2}} \Lambda_{\Z_{1/2}}[e^{(1)},\ldots,e^{(l)}]
		\end{align*}
		where $e^{(j)}, 1\leq j\leq l$ is the Euler class of the sphere $S^{2m_i}$, and can be identified as the Euler class of the $j$-th universal vector bundle over ${SO(n)}/\big({SO(2m_1)\times \cdots \times SO(2m_k)}\big)$ by Lemma\,\ref{lem:BundleTower}\,(b). To specify the $H_{T^M}^*(X,{\Z_{1/2}})$-algebra structure, we lift $e^{(j)}$ to be the equivariant Euler class $\tilde{e}^{(j)}$, and lift $(e^{(j)})^2=0$ to be $(\tilde{e}^{(j)})^2=\tilde{p}^{(j)}_{m_{j}}$ where $\tilde{p}^{(j)}_{m_{j}}$ is the equivariant top Pontryagin class of the $j$-th universal bundle over ${SO(n)}/\big({SO(2m_1)\times \cdots \times SO(2m_k)}\big)$. Hence, as an $H_{T^M}^*(X,{\Z_{1/2}})$-algebra,
		\begin{align*}
		 H_{T^M}^*\Big(\frac{SO(n)}{SO(2m_1)\times \cdots\times SO(2m_k)},{\Z_{1/2}}\Big)=\frac{H_{T^M}^*(X,{\Z_{1/2}})[\tilde{e}^{(j)}\mid_{1\leq j\leq l}]}{\big\langle(\tilde{e}^{(j)})^2=\tilde{p}^{(j)}_{m_{j}}\mid_{1\leq j\leq l}\big\rangle}.
		\end{align*}
		Comparing the above expression with Theorem\,\ref{thm:orientFlag1}\,(b) and using dimension counting, we then get the expression for $H_{T^M}^*(X,{\Z_{1/2}})$.
\end{proof}

\begin{thm}\label{thm:realFlag2}
	As a ${\Z_{1/2}}[\alpha_1,\ldots,\alpha_M]$-algebra, $H_{T^M}^*(Y,{\Z_{1/2}})$ is
		\begin{align*}
		\frac{{\Z_{1/2}}[\alpha_l\mid_{1\leq l\leq M}][\tilde{p}_i^{(j)}\mid_{1\leq j\leq k,1\leq i\leq m_j}]}{\big\langle\prod_{j=1}^{k}\tilde{p}^{(j)}=\prod_{l=1}^{M}(1+\alpha^2_l) \big\rangle} \otimes_{\Z_{1/2}} H^*\Big( \frac{SO(n)}{SO(2M+1)},{\Z_{1/2}}\Big).
		\end{align*}
\end{thm}
\begin{proof}
	Because of the $S\big(O(1)^k\big)$-covering: $X\rightarrow Y$, we would try to derive $H_{T^M}^*(Y,{\Z_{1/2}})$ as the $S\big(O(1)^k\big)$-invariant subring of $H_{T^M}^*(X,{\Z_{1/2}})$. However, it is not immediate to know the $S\big(O(1)^k\big)$-action on the exterior algebra $H^*\big( {SO(n)}/{SO(2M+1)},{\Z_{1/2}}\big)$. Instead, we try the following steps.
		\begin{enumerate}[(i)]
			\item Consider the following $T^M$-equivariant bundle projection:
			\[
			\frac{O(n)}{\big(O(2m_1)\times O(1)\big) \times\cdots\times \big(O(2m_l)\times O(1)\big)\times O(2m_{l+1})\times\cdots\times O(2m_k)} \longrightarrow Y
			\]
			whose fibre is 
			\[
			\frac{O(2m_1+1)}{O(2m_1)\times O(1)}\times\cdots\times \frac{O(2m_l+1)}{O(2m_l)\times O(1)} \cong \R P^{2m_1}  \times\cdots\times \R P^{2m_l}
			\]
			which has trivial cohomology in $\Z_{1/2}$ coefficients. Then, as equivariant cohomology ring,
			$$
			H_{T^M}^*(Y,{\Z_{1/2}})=H_{T^M}^*\Big(\frac{O(n)}{O(2m_1)\times\cdots\times O(2m_k)\times O(1)^l},{\Z_{1/2}}\Big)
			$$
			where we have rearranged the positions of the $O(1)^l$.
			\item Consider the $T^M$-equivariant covering map
			\begin{align*}
			\frac{O(n)}{O(2m_1)\times\cdots\times O(2m_k)} \longrightarrow \frac{O(n)}{O(2m_1)\times\cdots\times O(2m_k)\times O(1)^l}
			\end{align*} 
			and note that the covering group $O(1)^l$ operates in different blocks from $O(2m_1)\times\cdots\times O(2m_k)$. Since $S\big(O(1)^l\big) \subset SO(l)$, every element of $S\big(O(1)^l\big)$ is joined to the identity via a path in $SO(l)$, hence is homotopic to the identity as homomorphisms on $H_{T^M}^*\big({O(n)}/\big(O(2m_1)\times\cdots\times O(2m_k)\big),{\Z_{1/2}}\big)$. Then the $O(1)^l$ action on $H_{T^M}^*\big({O(n)}/\big(O(2m_1)\times\cdots\times O(2m_k)\big),{\Z_{1/2}}\big)$ is concentrated in the action of a single component $O(1)\cong O(1)^l/S\big(O(1)^l\big)$. Using Lemma\,\ref{lem:coveringSpace} twice, we get the ring identifications
			\begin{align*}
			H_{T^M}^*\Big(\frac{O(n)}{O(2m_1)\times\cdots\times O(2m_k)\times O(1)^l},{\Z_{1/2}}\Big)
			&= H_{T^M}^*\Big(\frac{O(n)}{O(2m_1)\times\cdots\times O(2m_k)},{\Z_{1/2}}\Big)^{O(1)^l}\\
			&= H_{T^M}^*\Big(\frac{O(n)}{O(2m_1)\times\cdots\times O(2m_k)},{\Z_{1/2}}\Big)^{O(1)}\\
			&= H_{T^M}^*\Big(\frac{O(n)}{O(2m_1)\times\cdots\times O(2m_k)\times O(1)},{\Z_{1/2}}\Big).
			\end{align*}
			\item Consider the $T^M$-equivariant covering map
			\[
			\frac{SO(n)}{SO(2m_1)\times \cdots \times SO(2m_k)\times SO(1)} \longrightarrow \frac{SO(n)}{S\Big(O(2m_1)\times \cdots \times O(2m_k)\times O(1)\Big)}
			\]
			with the covering group $S\big(O(1)^{k+1}\big)=\big\{(\epsilon_1,\ldots,\epsilon_{k+1})\mid \epsilon_j=\pm 1,\, \prod_{j=1}^{k+1}\epsilon_j=1\big\}\cong \Z_2^k$. Using Lemma\,\ref{lem:coveringSpace}, we have the ring identification
			\[
			H_{T^M}^*\Big(\frac{O(n)}{O(2m_1)\times\cdots\times O(2m_k)\times O(1)},{\Z_{1/2}}\Big)=
			H_{T^M}^*\Big(\frac{SO(n)}{SO(2m_1)\times \cdots \times SO(2m_k)},{\Z_{1/2}}\Big)^{S\big(O(1)^{k+1}\big)}.
			\]
			Similar to the argument in Theorem\,\ref{thm:realFlag1}\,(b), for any $\tilde{e}^{(j)}\in H_{T^M}^*\big({SO(n)}/{SO(2m_1)\times \cdots \times SO(2m_k)},{\Z_{1/2}}\big)$, where $1\leq j\leq k$, set $\epsilon_j=\epsilon_{k+1}=-1$ while all the other $\epsilon_j=1$, then 
			\[
			(\epsilon_1,\ldots,\epsilon_{k+1}) \Cdot \tilde{e}^{(j)}=\epsilon_j \tilde{e}^{(j)} = -\tilde{e}^{(j)}.
			\]
			Hence, any product of odd number of equivariant Euler classes is not $S\big(O(1)^{k+1}\big)$-invariant. In contrast, the Pontryagin classes are $S\big(O(1)^{k+1}\big)$-invariant. Since $S\big(O(1)^{k+1}\big)\subset S\big(O(2m_1)\times \cdots \times O(2m_k)\times O(1)\big) \subset SO(2M+1)$, the $S\big(O(1)^{k+1}\big)$-action on ${SO(n)}/{SO(2M+1)}$ is trivial, then the exterior algebra $H^*({SO(n)}/{SO(2M+1)},{\Z_{1/2}})$ is also $S\big(O(1)^{k+1}\big)$-invariant. 
		\end{enumerate}
\end{proof}

\begin{cor}\label{cor:flag2}
	The $\Z_{1/2}$-coefficient ordinary cohomology rings of the oriented and real flag manifolds $X,Y$ in Theorems\,\ref{thm:orientFlag2},\,\ref{thm:realFlag2} can be obtained by substituting $\alpha_1=\cdots=\alpha_M=0$. The ordinary Poincare polynomials of $X$ and $Y$ are respectively
	\[
	\binom{M}{m_1,\ldots,m_k}_{t^4}\cdot\prod_{j=l+1}^k (1+t^{2m_j})\cdot P\Big( \frac{SO(n)}{SO(2M+1)},t\Big), \qquad \quad
	\binom{M}{m_1,\ldots,m_k}_{t^4} \cdot P\Big( \frac{SO(n)}{SO(2M+1)},t\Big).
	\]
	As a result, if $0<l< k$, i.e. $(n_1,\ldots,n_k)$ is neither all even nor all odd, then $Y$ is non-orientable. 
\end{cor}
\begin{proof}
	The statement about ordinary cohomology and Poincare polynomials is similar to Corollary\,\ref{cor:flag1}. For the orientability, $Y$ is non-orientable $\Leftrightarrow$ the above two polynomials have different degrees $\Leftrightarrow$ the term $\prod_{j=l+1}^k (1+t^{2m_j})$ is nontrivial $\Leftrightarrow$ $0<l<k$ $\Leftrightarrow$ $(n_1,\ldots,n_k)$ is neither all even nor all odd. 
\end{proof}

\begin{exm}
	Consider the real or oriented Grassmannians
	{
		\begin{align*}
		G_k(\R^n)=Fl(k,n-k,\R^n)=\frac{O(n)}{O(k)\times O(n-k)},     \qquad  G^o_k(\R^n)=Fl^o(k,n-k,\R^n)=\frac{SO(n)}{SO(k)\times SO(n-k)}
		\end{align*}}among which only the even-dimensional Grassmannians
	{
		\begin{align*}
		G_{2k}(\R^{2n})&=\frac{O(2n)}{O(2k)\times O(2n-2k)}    &  G^o_{2k}(\R^{2n})&=\frac{SO(2n)}{SO(2k)\times SO(2n-2k)}\\
		G_{2k}(\R^{2n+1})&=\frac{O(2n+1)}{O(2k)\times O(2n-2k+1)}    &  G^o_{2k}(\R^{2n+1})&=\frac{SO(2n+1)}{SO(2k)\times SO(2n-2k+1)}\\
		G_{2k+1}(\R^{2n+1})&=\frac{O(2n+1)}{O(2k+1)\times O(2n-2k)}    &  G^o_{2k+1}(\R^{2n+1})&=\frac{SO(2n+1)}{SO(2k+1)\times SO(2n-2k)}
		\end{align*}}satisfy the equal-rank condition. Noticing the identifications $G_{2k}(\R^{2n+1})\cong G_{2n-2k+1}(\R^{2n+1}),\,G^o_{2k}(\R^{2n+1})\cong G^o_{2n-2k+1}(\R^{2n+1})$ between the second the third rows, we only consider the first two rows of Grassmannians and get the equivariant cohomology rings in $\Z_{1/2}$ coefficients (see\,\cite{GHV76} p.\,478,\,480 for the ordinary cohomology of the oriented Grassmannians and\,\cite{CK,Sa17,Carl,He} for more details):
	{
		\begin{gather*}
		H_{T^n}^*(G_{2k}(\R^{2n}),\Z_{1/2})\cong H_{T^n}^*(G_{2k}(\R^{2n+1}),\Z_{1/2})=\frac{\Z_{1/2}[\alpha_1,\dots,\alpha_n;\tilde{p}_1,\ldots,\tilde{p}_k;\tilde{p}'_1,\ldots,\tilde{p}'_{n-k}]}{\big\langle \tilde{p}\tilde{p}'=\prod_{i=1}^{n}(1+\alpha_i^2)\big\rangle}\\
		H_{T^n}^*(G^o_{2k}(\R^{2n}),\Z_{1/2})=\frac{\Z_{1/2}[\alpha_1,\dots,\alpha_n;\tilde{p}_1,\ldots,\tilde{p}_k;\tilde{p}'_1,\ldots,\tilde{p}'_{n-k};\tilde{e},\tilde{e}']}{\big\langle \tilde{p}\tilde{p}'=\prod_{i=1}^{n}(1+\alpha_i^2);\tilde{e}^2=\tilde{p}_k,\tilde{e}'^2=\tilde{p}'_{n-k};\tilde{e}\tilde{e}'=\prod_{i=1}^{n}\alpha_i\big\rangle} \\
		H_{T^n}^*(G^o_{2k}(\R^{2n+1}),\Z_{1/2})=\frac{\Z_{1/2}[\alpha_1,\dots,\alpha_n;\tilde{p}_1,\ldots,\tilde{p}_k;\tilde{p}'_1,\ldots,\tilde{p}'_{n-k};\tilde{e}]}{\big\langle \tilde{p}\tilde{p}'=\prod_{i=1}^{n}(1+\alpha_i^2);\tilde{e}^2=\tilde{p}_k\big\rangle}
		\end{gather*}}where $\tilde{p}=1+\tilde{p}_1+\cdots \tilde{p}_k,\,\tilde{p}'=1+\tilde{p}'_1+\cdots \tilde{p}'_{n-k},\tilde{e},\tilde{e}'$ are the equivariant total Pontryagin and the equivariant Euler classes of the universal bundle and its complementary bundle over those even-dimensional Grassmannians. The Poincare polynomials of the even-dimensional real or oriented Grassmannians are
	{
		\begin{align*}
		P(G_{2k}(\R^{2n}),t)=P(G_{2k}(\R^{2n+1}),t)=\binom{n}{k,n-k}_{t^4},\\
		P(G^o_{2k}(\R^{2n}),t)=\binom{n}{k,n-k}_{t^4}\cdot\frac{(1+t^{2k})(1+t^{2n-2k})}{1+t^{2n}}, \qquad  P(G^o_{2k}(\R^{2n+1}),t)&=\binom{n}{k,n-k}_{t^4}\cdot(1+t^{2k}).
		\end{align*}}
\end{exm}

\begin{exm}
	For the odd-dimensional Grassmannians 
	{
		\begin{align*}
		G_{2k+1}(\R^{2n+2})=\frac{O(2n+2)}{O(2k+1)\times O(2n-2k+1)}    \qquad \qquad G^o_{2k+1}(\R^{2n+2})=\frac{SO(2n+2)}{SO(2k+1)\times SO(2n-2k+1)},
		\end{align*}}the denominators have the maximal torus $T^n=SO(2)^k \times SO(2)^{n-k}$ of rank less than the rank of the numerators. The $\Z_{1/2}$-coefficient equivariant cohomology rings of both Grassmannians are similar to the case of even-dimensional Grassmannians except for an exterior generator (see\,\cite{Ta62},\,\cite{GHV76} p.\,480,\,\cite{CK} for the ordinary cohomology, and\,\cite{Sa17,Carl,He} for the equivariant cohomology): 
	{
		\begin{align*}
		\frac{{\Z_{1/2}}[\alpha_1,\dots,\alpha_n;\tilde{p}_1,\ldots,\tilde{p}_k;\tilde{p}'_1,\ldots,\tilde{p}'_{n-k}]}{\big\langle \tilde{p}\tilde{p}'=\prod_{i=1}^{n}(1+\alpha_i^2)\big\rangle}\otimes_{\Z_{1/2}} \Lambda_{\Z_{1/2}}[x_{n+1}]
		\end{align*}}where $x_{n+1}$ is of degree $2n+1$. The Poincare polynomials of both odd-dimensional Grassmannians are
	{
		\begin{align*}
		\binom{n}{k,n-k}_{t^4}\cdot(1+t^{2n+1}).
		\end{align*}}
\end{exm}

\begin{exm}
	Neither of the real and oriented complete flag manifolds
	{
		\begin{align*}
		Fl(\R^n)=\frac{O(n)}{O(1)^n}=\frac{O(n)}{\Z_2^n}  \qquad  \qquad\qquad \qquad Fl^o(\R^n)=\frac{SO(n)}{SO(1)^n}=SO(n)
		\end{align*}}satisfy the equal-rank condition. But we can use the Step\,(ii) of Theorem\,\ref{thm:realFlag2} to see that 
	$$H^*\Big(\frac{O(n)}{O(1)^n},\Z_{1/2}\Big)= H^*\Big(\frac{O(n)}{O(1)},\Z_{1/2}\Big)=H^*(SO(n),\Z_{1/2})$$ which is an exterior algebra as in Theorem\,\ref{thm:realStief}. This identification was kindly communicated by Liviu Mare to the author and largely inspired the Step\,(ii) of Theorem\,\ref{thm:realFlag2}. 
\end{exm}
\begin{exm}
	Consider the following equal-rank flag manifolds
	{
		\begin{align*}
		&Fl(\underbrace{2,\ldots,2}_{n \text{ items}},\R^{2n})=\frac{O(2n)}{O(2)^n}    &  &Fl^o(\underbrace{2,\ldots,2}_{n \text{ items}},\R^{2n})=\frac{SO(2n)}{SO(2)^n}\\
		&Fl(\underbrace{2,\ldots,2}_{n \text{ items}},1,\R^{2n+1})=\frac{O(2n+1)}{O(2)^n \times O(1)}   &  &Fl^o(\underbrace{2,\ldots,2}_{n \text{ items}},1,\R^{2n+1})=\frac{SO(2n+1)}{SO(2)^n \times SO(1)}
		\end{align*} }on which there are the canonical left actions of $T^n=SO(2)^n$. We have the equivariant cohomology rings:
	{
		\begin{align*}
		H_{T^n}^*\Big(\frac{O(2n)}{O(2)^n},{\Z_{1/2}}\Big)&= \frac{{\Z_{1/2}}[\alpha_1,\dots,\alpha_n;t^2_1,\ldots,t^2_n]}{\big\langle\prod_{i=1}^{n}(1+t^2_i)=\prod_{i=1}^{n}(1+\alpha^2_i)\big\rangle}\\
		H_{T^n}^*\Big(\frac{SO(2n)}{SO(2)^n},{\Z_{1/2}}\Big)&=\frac{{\Z_{1/2}}[\alpha_1,\dots,\alpha_n;t_1,\ldots,t_n]}{\big\langle\prod_{i=1}^{n}(1+t^2_i)=\prod_{i=1}^{n}(1+\alpha^2_i);\prod_{i=1}^{n}t_i=\prod_{i=1}^{n}\alpha_i\big\rangle}\\
		H_{T^n}^*\Big(\frac{O(2n+1)}{O(2)^n\times O(1)},{\Z_{1/2}}\Big)&=\frac{{\Z_{1/2}}[\alpha_1,\dots,\alpha_n;t^2_1,\ldots,t^2_n]}{\big\langle\prod_{i=1}^{n}(1+t^2_i)=\prod_{i=1}^{n}(1+\alpha^2_i)\big\rangle}\\
		H_{T^n}^*\Big(\frac{SO(2n+1)}{SO(2)^n },{\Z_{1/2}}\Big)&=\frac{{\Z_{1/2}}[\alpha_1,\dots,\alpha_n;t_1,\ldots,t_n]}{\big\langle\prod_{i=1}^{n}(1+t^2_i)=\prod_{i=1}^{n}(1+\alpha^2_i)\big\rangle}
		\end{align*}}where for the real flag manifolds $O(2n)/O(2)^n,O(2n+1)/O(2)^n \times O(1)$, the $t^2_i$ is the equivariant first Pontryagin class of the $i$-th universal $\R^2$-bundle; and for the oriented flag manifolds $SO(2n)/SO(2)^n,SO(2n+1)/SO(2)^n$, the $t_i$ is the equivariant Euler class of the $i$-th universal oriented $\R^2$-bundle. The Poincare polynomials are
	{
		\begin{align*}
		P\Big(\frac{O(2n)}{O(2)^n},t\Big)=P\Big(\frac{O(2n+1)}{O(2)^n\times O(1)},t\Big)=[n]_{t^4}!,\\
		P\Big(\frac{SO(2n)}{SO(2)^n},t\Big)=[n]_{t^4}!\cdot \frac{(1+t^2)^n}{1+t^{2n}},\qquad \qquad
		P\Big(\frac{SO(2n+1)}{SO(2)^n},t\Big)&=[n]_{t^4}!\cdot(1+t^2)^n.
		\end{align*}}
\end{exm}

\vskip 20pt
\section{Integral equivariant cohomology rings of the complex or quaternionic flag manifolds}
\vskip 15pt
The complex or quaternionic counterparts of the real or oriented flag manifolds can be defined as 
\[
Fl(n_1,\ldots,n_k,\C^n)\triangleq \frac{U(n)}{U(n_1)\times\cdots\times U(n_k)},
\qquad \qquad  Fl(n_1,\ldots,n_k,\mathbb{H}^n)\triangleq\frac{Sp(n)}{Sp(n_1)\times\cdots\times Sp(n_k)}
\]
where $n\geq n_1+\cdots+n_k$. These two cases are much more pleasant than the real or oriented cases because $U(n),\,Sp(n)$ do not have torsions in their integral cohomology, nor will the parities of $n_1,\ldots,n_k,n$ cause any losses of the ranks of the maximal tori.

Write $N= n_1+\cdots+n_k$ and set $G=U,\,Sp$ and $\D=\C,\,\mathbb{H}$ accordingly, and consider the $T^N$-equivariant bundle that we have seen from the real and oriented cases where the fibre is a usual flag manifold and the base is a Stiefel manifold.
\begin{figure}[H]
	\centering
	\begin{minipage}[b]{0.45\textwidth}
		\centering
		\begin{tikzpicture}[description/.style={fill=white,inner sep=2pt}]
		\matrix (m) [matrix of math nodes, row sep=2em,
		column sep=2em, text height=1.5ex, text depth=0.25ex]
		{ \frac{G(N)}{G(n_1)\times \cdots \times G(n_k)} & \frac{G(n)}{G(n_1)\times \cdots \times G(n_k)} \\
			& \frac{G(n)}{G(N)}\\};
		\path[right hook->] 
		(m-1-1) edge node[above] {$\iota$} (m-1-2);
		\path[->]
		(m-1-2) edge node[left] {$\pi$} (m-2-2);
		\end{tikzpicture}
	\end{minipage}
	\begin{minipage}[b]{0.45\textwidth}
		\centering
		\begin{tikzpicture}[description/.style={fill=white,inner sep=2pt}]
		\matrix (m) [matrix of math nodes, row sep=2em,
		column sep=2em, text height=1.5ex, text depth=0.25ex]
		{ Fl(n_1,\ldots,n_k,\D^N) & Fl(n_1,\ldots,n_k,\D^n) \\
			& V_{n-N}(\D^n).\\};
		\path[right hook->] 
		(m-1-1) edge node[above] {$\iota$} (m-1-2);
		\path[->]
		(m-1-2) edge node[left] {$\pi$} (m-2-2);
		\end{tikzpicture}
	\end{minipage}
\end{figure}

Since the fibre satisfies the equal-rank condition and the involving groups do not have torsions in their integral cohomology, we use the equivariant version of the Leray-Borel Theorem\,\ref{thm:BundleLerayBorel}, the complex or quaternionic versions of Theorems\,\ref{thm:realStief},\,\ref{thm:equivRealStiefel} and Relations\,\eqref{eq:EquivWhitney},\,\eqref{eq:EquivsqEuler},\,\eqref{eq:EquivWhitneyEuler} to get

\begin{thm}\label{thm:equivCplxFlag}
	The $H^*_{T^N}(pt,\Z)=\Z[\alpha_1,\ldots,\alpha_{N}]$-algebra structures of $H^*_{T^N}({G(N)}/{G(n_1)\times \cdots \times G(n_k)},\Z)$ for the complex or quaternionic flag manifolds respectively are:
	\begin{align*}
	\frac{{\Z}[\alpha_l\mid_{1\leq l\leq N}][\tilde{c}_i^{(j)}\mid_{1\leq j\leq k,1\leq i\leq n_j}]}{\big\langle\prod_{j=1}^{k}\tilde{c}^{(j)}=\prod_{l=1}^{N}(1+\alpha_l) \big\rangle} \otimes_{\Z} \Lambda_{\Z}[x_m \mid_{N+1\leq m\leq n}],\quad 
	\frac{{\Z}[\alpha_l\mid_{1\leq l\leq N}][\tilde{q}_i^{(j)}\mid_{1\leq j\leq k,1\leq i\leq n_j}]}{\big\langle\prod_{j=1}^{k}\tilde{q}^{(j)}=\prod_{l=1}^{N}(1+\alpha^2_l) \big\rangle} \otimes_{\Z} \Lambda_{\Z}[y_m \mid_{N+1\leq m\leq n}]
	\end{align*}
	where $\tilde{c}^{(j)}=1+\tilde{c}_1^{(j)}+\cdots+\tilde{c}_{n_j}^{(j)}$ and $q^{(j)}=1+\tilde{q}_1^{(j)}+\cdots+\tilde{q}_{n_j}^{(j)}$ are the equivariant total Chern and quaternionic classes of the $j$-th universal vector bundles, and $x_m,y_m$ are of degrees $2m-1,4m-1$. The ordinary cohomology rings can be obtained by substituting $\alpha_1=\cdots=\alpha_k=0$, and the Poincare polynomials are respectively
	\[
	\binom{N}{n_1,\ldots,n_k}_{t^2} \cdot \prod_{i=N+1}^{n}(1+t^{2i-1}),\qquad \qquad \qquad  \binom{N}{n_1,\ldots,n_k}_{t^4} \cdot \prod_{i=N+1}^{n}(1+t^{4i-1}).
	\]
\end{thm}

\begin{rmk}
	The Leray-Borel descriptions of the ordinary cohomology rings of the complex or quaternionic flag manifolds have been explicitly given by Borel\,\cite{Bo53a} (also see \cite{GHV76} Chap.\,XI). The Leray-Borel descriptions of the equivariant cohomology rings of complex Grassmannians and complete flag manifolds have been given by Goldin\,\cite{Go99}, and the case of quaternionic complete flag manifold has been given by Mare\,\cite{Ma08}.
\end{rmk}

\subsection{Connections with the Schubert calculus}
In the following discussion, assume $n=n_1+\cdots+n_k$. The cohomology ring of the usual complex flag manifold $Fl(n_1,\ldots,n_k,\C^n)$ is generated by the Chern classes and does not contain the exterior part. On the other hand, for a complex semisimple Lie group $G$ and its parabolic subgroup $P$, the Bruhat-Chevalley decomposition gives an additive basis consisting of the Schubert classes $\big\{\mathcal{S}_\lambda\mid \lambda \in W_G/W_P\big\}$ for the cohomology of the complex flag variety $G/P$. The ring structure of the cohomology can be described as $\mathcal{S}_\lambda \cdot \mathcal{S}_\mu = \sum_{\nu} C_{\lambda,\mu}^\nu \mathcal{S}_\nu $
where $C_{\lambda,\mu}^\nu \in \Z$ are the Littlewood-Richardson coefficients (see\,\cite{Co} for a geometric rule of the LR coefficients of the complex flag manifolds). 

The work of Bernstein, Gel'fand and Gel'fand\,\cite{BGG73} and Demazure\,\cite{De74} establishes the connections between these two descriptions which are also valid in the equivariant setting as shown by Arabia\,\cite{Ar89}. Hence, for the complex flag manifold $Fl(n_1,\ldots,n_k,\C^n)=U(n_1+\cdots+n_k)/\big(U(n_1)\times\cdots\times U(n_k)\big)$, there exists a ring identification between the two descriptions of its equivariant cohomology ring, one in equivariant Chern classes and the other one in equivariant Schubert classes $\tilde{\mathcal{S}}_\lambda$ indexed by $\lambda \in W_G/W_P=S_{n}/(S_{n_1}\times\cdots\times S_{n_k})$ where $S_{n}$ is the symmetric group of $n$ elements:
\begin{align*}
\mathcal{A}:\frac{{\Z}[\alpha_l\mid_{1\leq l\leq n}][\tilde{c}_i^{(j)}\mid_{1\leq j\leq k,1\leq i\leq n_j}]}{\big\langle\prod_{j=1}^{k}\tilde{c}^{(j)}=\prod_{l=1}^{n}(1+\alpha_l) \big\rangle}\overset{=}{\longrightarrow} \frac{\Z[\alpha_l\mid_{1\leq l\leq n}][\tilde{\mathcal{S}}_\lambda\mid_{\lambda \in S_{n}/(S_{n_1}\times\cdots\times S_{n_k})}]}{\big\langle\tilde{\mathcal{S}}_\lambda \cdot \tilde{\mathcal{S}}_\mu = \sum_{\nu} \tilde{C}_{\lambda,\mu}^\nu \tilde{\mathcal{S}}_\nu\big\rangle}
\end{align*}
where $\tilde{C}_{\lambda,\mu}^\nu\in \Z[\alpha_1,\ldots,\alpha_n]$ are the equivariant LR coefficients. Kaji\,\cite{Ka} gives an explicit identification of these two descriptions in the case of complex complete flag manifold.

For the quaternionic flag manifold $Fl(n_1,\ldots,n_k,\mathbb{H}^n)=Sp(n_1+\cdots+n_k)/\big(Sp(n_1)\times\cdots\times Sp(n_k)\big)$, its equivariant cohomology ring can be viewed a ${\Z}[\alpha_1,\ldots,\alpha_k]$-subalgebra of $H^*_{T^n}\big(Fl(n_1,\ldots,n_k,\mathbb{C}^n),\Z\big)$, such that for each $1\leq j \leq k$, the alternating sum $1+\sum_{i=1}^{n_j}(-1)^i\tilde{q}_i^{(j)}$ is identified with $\big(1+\sum_{i=1}^{n_j}\tilde{c}_i^{(j)}\big)\big(1+\sum_{i=1}^{n_j}(-1)^i\tilde{c}_i^{(j)}\big)$. On the other hand, $H^*_{T^n}\big(Fl(n_1,\ldots,n_k,\mathbb{H}^n),\Z\big)$ is the result of a degree-doubling operation on $H^*_{T^n}\big(Fl(n_1,\ldots,n_k,\mathbb{C}^n),\Z\big)$:
\begin{align*}
\mathcal{D}:\frac{{\Z}[\alpha_l\mid_{1\leq l\leq n}][\tilde{c}_i^{(j)}\mid_{1\leq j\leq k,1\leq i\leq n_j}]}{\big\langle\prod_{j=1}^{k}\tilde{c}^{(j)}=\prod_{l=1}^{n}(1+\alpha_l) \big\rangle}\longrightarrow \frac{{\Z}[\alpha_l\mid_{1\leq l\leq N}][\tilde{q}_i^{(j)}\mid_{1\leq j\leq k,1\leq i\leq n_j}]}{\big\langle\prod_{j=1}^{k}\tilde{q}^{(j)}=\prod_{l=1}^{N}(1+\alpha^2_l) \big\rangle}: \tilde{c}_i^{(j)}\longmapsto \tilde{q}_i^{(j)},\, \alpha_l\longmapsto \alpha_l^2
\end{align*} 
which is a ring homomorphism, but not a ${\Z}[\alpha_1,\ldots,\alpha_k]$-algebra homomorphism. Applying the operation $\mathcal{D}$ on the two descriptions of $H^*_{T^n}\big(Fl(n_1,\ldots,n_k,\mathbb{C}^n),\Z\big)$, we get a ring identification between the two descriptions of $H^*_{T^n}\big(Fl(n_1,\ldots,n_k,\mathbb{H}^n),\Z\big)$:
\begin{align*}
\mathcal{A}_\mathcal{D}:\frac{{\Z}[\alpha_l\mid_{1\leq l\leq n}][\tilde{q}_i^{(j)}\mid_{1\leq j\leq k,1\leq i\leq n_j}]}{\big\langle\prod_{j=1}^{k}\tilde{q}^{(j)}=\prod_{l=1}^{n}(1+\alpha^2_l) \big\rangle}\overset{=}{\longrightarrow}\frac{\Z[\alpha_l\mid_{1\leq l\leq n}][\mathcal{D}(\tilde{\mathcal{S}}_\lambda)\mid_{\lambda \in S_{n}/(S_{n_1}\times\cdots\times S_{n_k})}]}{\big\langle\mathcal{D}(\tilde{\mathcal{S}}_\lambda) \cdot \mathcal{D}(\tilde{\mathcal{S}}_\mu) = \sum_{\nu} \mathcal{D}(\tilde{C}_{\lambda,\mu}^\nu) \mathcal{D}(\tilde{\mathcal{S}}_\nu)\big\rangle}
\end{align*}
where the equivariant LR coefficients $\mathcal{D}(\tilde{C}_{\lambda,\mu}^\nu)$ are obtained from the polynomials $\tilde{C}_{\lambda,\mu}^\nu$ by substituting $\alpha_l^2$ for $\alpha_l$, and the ordinary LR coefficients $\mathcal{D}({C}_{\lambda,\mu}^\nu)$ are the same as ${C}_{\lambda,\mu}^\nu$.

Similarly, for the real flag manifold $Fl(2n_1,\ldots,2n_k,\mathbb{R}^{2n})=O(2n_1+\cdots+2n_k)/\big(O(2n_1)\times\cdots\times O(2n_k)\big)$, we can work in $\Z_{1/2}$ coefficients and replace the above equivariant quaternionic class $\tilde{q}_i^{(j)}$ by the equivariant Pontryagin class $\tilde{p}_i^{(j)}$ and get the same equivariant LR coefficient $\mathcal{D}(\tilde{C}_{\lambda,\mu}^\nu)$. However, for oriented flag manifolds and other types of real flag manifolds where we don't have all the even integers $2n_1,\ldots,2n_k$, there need to associate Schubert symbols to the equivariant Euler classes and the exterior generators. 

For the real flag manifold $Fl(n_1,\ldots,n_k,\mathbb{R}^{n})=O(n_1+\cdots+n_k)/\big(O(n_1)\times\cdots\times O(n_k)\big)$, its mod-$2$ $\Z_2^n$-equivariant cohomology ring is the result of a degree-halving operation on $H^*_{T^n}\big(Fl(n_1,\ldots,n_k,\mathbb{C}^n),\F_2\big)$:
\begin{align*}
\mathcal{H}:\frac{{\F_2}[\alpha_l\mid_{1\leq l\leq n}][\tilde{c}_i^{(j)}\mid_{1\leq j\leq k,1\leq i\leq n_j}]}{\big\langle\prod_{j=1}^{k}\tilde{c}^{(j)}=\prod_{l=1}^{n}(1+\alpha_l) \big\rangle}\longrightarrow \frac{{\F_2}[\beta_l\mid_{1\leq l\leq N}][\tilde{w}_i^{(j)}\mid_{1\leq j\leq k,1\leq i\leq n_j}]}{\big\langle\prod_{j=1}^{k}\tilde{w}^{(j)}=\prod_{l=1}^{N}(1+\beta_l) \big\rangle}: \tilde{c}_i^{(j)}\longmapsto \tilde{w}_i^{(j)},\, \alpha_l\longmapsto \beta_l.
\end{align*} 
Applying the operation $\mathcal{H}$ on the two descriptions of $H^*_{T^n}\big(Fl(n_1,\ldots,n_k,\mathbb{C}^n),\F_2\big)$, we get a ring identification between the two descriptions of $H^*_{\Z_2^n}\big(Fl(n_1,\ldots,n_k,\mathbb{R}^n),\F_2\big)$:
\begin{align*}
\mathcal{A}_\mathcal{H}:\frac{{\F_2}[\beta_l\mid_{1\leq l\leq n}][\tilde{w}_i^{(j)}\mid_{1\leq j\leq k,1\leq i\leq n_j}]}{\big\langle\prod_{j=1}^{k}\tilde{w}^{(j)}=\prod_{l=1}^{n}(1+\beta_l) \big\rangle}\overset{=}{\longrightarrow}\frac{\F_2[\beta_l\mid_{1\leq l\leq n}][\mathcal{H}(\tilde{\mathcal{S}}_\lambda)\mid_{\lambda \in S_{n}/(S_{n_1}\times\cdots\times S_{n_k})}]}{\big\langle\mathcal{H}(\tilde{\mathcal{S}}_\lambda) \cdot \mathcal{H}(\tilde{\mathcal{S}}_\mu) = \sum_{\nu} \mathcal{H}(\tilde{C}_{\lambda,\mu}^\nu) \mathcal{H}(\tilde{\mathcal{S}}_\nu)\big\rangle}
\end{align*}
where the mod-$2$ equivariant LR coefficients $\mathcal{H}(\tilde{C}_{\lambda,\mu}^\nu)$ are obtained from the polynomial $\tilde{C}_{\lambda,\mu}^\nu$ by substituting $\beta_l$ for $\alpha_l$, and the mod-$2$ ordinary LR coefficients $\mathcal{H}({C}_{\lambda,\mu}^\nu)$ are the same as ${C}_{\lambda,\mu}^\nu$.

\vskip 20pt
\bibliographystyle{amsalpha}

\end{document}